%% file: NGA5.TEX
\documentclass[a4paper,11pt,english]{article}

\usepackage{babel}
\usepackage{amsfonts}
\usepackage{amssymb}
\usepackage{oldlfont} 
\usepackage{tocloft}

\input{tcilatex}
\input{commands2e}
\renewcommand{\text}{\mbox}
\usepackage[colorlinks,
						bookmarks,
						pdftitle={Generalized Functions in Infinite Dimensional Analysis},
            pdfsubject={Generalized Functions in Infinite Dimensional Analysis},
            pdfkeywords={White Noise Analysis, Gaussian Analysis,
            Non-Gaussian Analysis, infinite dimensional Analysis},
            pdfauthor={Kondratiev, Streit, Westerkamp, Yan},
            plainpages=false,
            linkcolor={blue},
            citecolor={blue},
            urlcolor={blue}  		
						]{hyperref}

\newtheorem{theorem}{\sc Theorem}
\newtheorem{lemma}[theorem]{\sc Lemma}
\newtheorem{proposition}[theorem]{\sc Proposition}

\newtheorem{definition}[theorem]{\sc Definition}
\newtheorem{corollary}[theorem]{\sc Corollary}
\renewcommand{\example}{\noindent{\sc Example \stepcounter{example}%
                                     \arabic{example} }}
\renewcommand{\assumption}{\noindent{\sc Assumption }} 
\renewcommand{\proof}{\noindent{\sc Proof. \quad }} 
\renewcommand{\note}{\noindent{\sc Note. \quad }} 
\renewcommand{\remark}{\noindent{\sc Remark. \quad }} 
\renewcommand{\remarks}{\noindent{\sc Remarks. \quad }} 

\begin{document}

\author{Yuri G. {\sc Kondratiev} \smallskip \\Ludwig {\sc Streit} \smallskip %
\\Werner {\sc Westerkamp} \smallskip \\Jia-an {\sc Yan} }
\title{Generalized Functions in Infinite Dimensional Analysis}
\date{October 1995}
\maketitle

\begin{abstract}
We give a general approach to infinite dimensional non-Gaussian Analysis
which generalizes the work \cite{ADKS94} to measures which possess more
singular logarithmic derivative. This framework also includes the
possibility to handle measures of Poisson type.
\end{abstract}


\cftbeforesecskip=3pt
\tableofcontents

\section{Background and Introduction}

\input{n1intro}

\section{Preliminaries\label{Preliminaries}}

\subsection{Some facts on nuclear triples}

\input{nuclear}

\subsection{Holomorphy on locally convex spaces\label{Holomorphy}}

\input{holo}

\input{measures}\input{alms}\input{wick}\input{n1pos}\input{change}

\bigskip\ 

\noindent {\bf Acknowledgements. } 
The following results were presented at the workshop
``Mathematical Approach to Fluctuations" organized at IIAS 
by Professor T.~Hida. L.S.~would like to express his
gratitude to IIAS not only for the opportunity to do so, but also for
its academic hospitality in the weeks after the workshop which provided
a unique opportunity for continued work and collaboration.
The first author was supported by the
European program `Human capital and mobility' and by the `Deutsche
Forschungsgemeinschaft'. 
W.W.~gratefully acknowledges financial support of a scholarship from
`Graduiertenf\"orderung des Landes Nordrhein--Westfalen'. The last author
(J.A.Y) was partially supported by STRIDE and the National Natural Science
Foundation of China.

\TeXButton{TeX field}{\small \frenchspacing}

\medskip \normalsize
\hspace*{\fill} \begin{tabular}[t]{l}
                    Yuri G. Kondratiev$^{1,2}$ \smallskip 
                    \\ Ludwig Streit$^{1,3}$ \smallskip 
                    \\ Werner Westerkamp$^1$ \smallskip 
                    \\ Jia-an Yan$^{4}$ \\ \medskip\  \\ 
                  \end{tabular}   
      {\it        \begin{tabular}[t]{l}
                   $^1$BiBoS, Universit\"at Bielefeld, \\
                              D 33615 Bielefeld, Germany
                   \\[2mm] $^2$Institute of Mathematics, \\
                              Kiev, Ukraine
                   \\[2mm] $^3$CCM, Universidade da Madeira, \\
                              P9000 Funchal, Portugal
                   \\[2mm] $^4$Institute of Applied Mathematics, \\
                               Academica Sinica, \\
                               100080 Beijing, China
               \end{tabular}  }
\end{document}

%% file: commands2e.tex
\newcommand{\D}{\displaystyle}
\newcommand{\T}{\textstyle}
\newcounter{example}
\newcommand{\example}{\noindent{\bf Example \nopagebreak \stepcounter{example}%
                                     \arabic{example} }}
\renewcommand{\endproof}{\hspace*{\fill}$\Box$}

\newcommand{\C}{{\mathbb  C}}
\newcommand{\Ckl}{{\mathbb  C}} 
\newcommand{\Ckkl}{{\mathbb  C}}

\newcommand{\Z}{{\mathbb  Z}}
\newcommand{\R}{{\mathbb  R}}
\newcommand{\p}{{\mathbb  P}}
\newcommand{\Q}{{\mathbb  Q}}
\newcommand{\A}{{\mathbb  A}}
\newcommand{\N}{{\mathbb  N}}
\newcommand{\1}{{1\!\! 1}}
\newcommand{\E}{{\mathbb  E}}

\newcommand{\lnorm}{ \left| \! \left| \! \left| }
\newcommand{\rnorm}{ \right| \! \right| \! \right| }

\newcommand{\assumption}{\noindent{\bf Assumption }} 
\renewcommand{\note}{\noindent{\bf Note. \quad }} 
\newcommand{\remark}{\noindent{\bf Remark. \quad }} 
\newcommand{\remarks}{\noindent{\bf Remarks. \quad }} 

%% file: n1intro.tex

\noindent
White Noise Analysis and -- more generally -- Gaussian analysis have now
become of age, both date back approximately twenty years, for reviews we
refer to \cite{BeKo88,HKPS93}. Essential to both of them is an orthogonal
decomposition of the underlying $L^2$ space --- the ``chaos'' or ``Hermite''
or ``normal'' or ``multiple Wiener integral'' decomposition.\medskip\ 

One extension of this setup has been introduced by Y.M.~Berezansky: Starting
from certain field operators he constructs polynomial or orthogonal
decompositions with respect to the spectrum measures which need not
necessary be Gaussian, see e.g., \cite{BeLy93}.\medskip\ 

A different approach was recently proposed by \cite{AKS93}. For smooth
probability measures on infinite dimensional linear spaces a biorthogonal
decomposition is a natural extension of the orthogonal one that is well
known in Gaussian analysis. This biorthogonal ``Appell'' system has been
constructed for smooth measures by Yu.L.~Daletskii \cite{Da91}. For a
detailed description of its use in infinite dimensional analysis we refer to 
\cite{ADKS94}.\bigskip\ 

\noindent {\it Aim of the present work. }We consider the case of
non--degenerate measures on co-nuclear spaces with analytic characteristic
functionals. It is worth emphasizing that no further condition such as
quasi--invariance of the measure or smoothness of logarithmic derivatives
are required. The point here is that the important example of Poisson noise
is now accessible.

\noindent For any such measure $\mu $ we construct an Appell system $\A^\mu $
as a pair $(\p ^\mu ,\Q ^\mu )$ of Appell polynomials $\p^\mu $ and a
canonical system of generalized functions $\Q ^\mu $, properly associated to
the measure $\mu $.\bigskip\ 

\noindent {\it Central results. }Within the above framework

\begin{itemize}
\item  we obtain an explicit description of the test function space
introduced in \cite{ADKS94}

\item  in particular this space is in fact identical for all the measures
that we consider

\item  characterization theorems for generalized as well as test functions
are obtained analogously as in Gaussian analysis \cite{KLS94} for more
references see \cite{KLPSW94}

\item  the well known Wick product and the corresponding Wick calculus \cite
{KLS94} extends rather directly

\item  similarly, a full description of positive distributions (as measures)
will be given.
\end{itemize}

\noindent Finally we should like to underline here the important conceptual
role of holomorphy here as well as in earlier studies of Gaussian analysis
(see e.g., \cite{PS91,Ou91,KLPSW94,KLS94} as well as the references cited
therein).

%% file: nuclear.tex
\LaTeXparent{nga4.tex}

We start with a real separable Hilbert space ${\cal H}$ with inner product $%
(\cdot ,\cdot )$ and norm $\left| \cdot \right| $ . For a given separable
nuclear space ${\cal N}$ (in the sense of Grothendieck) densely
topologically embedded in ${\cal H}$ we can construct the nuclear triple{\tt %
\ } 
$$
{\cal N}\subset {\cal H\subset N^{\prime }}. 
$$
The dual pairing $\langle \cdot ,\cdot \rangle $ of ${\cal N}^{\prime }$ and 
${\cal N}$ then is realized as an extension of the inner product in ${\cal H}
$%
$$
\langle f,\,\xi \rangle =(f,\xi )\quad f\in {\cal H},\ \xi \in {\cal N} 
$$
Instead of reproducing the abstract definition of nuclear spaces (see e.g., 
\cite{Sch71}) we give a complete (and convenient) characterization in terms
of projective limits of Hilbert spaces.

\begin{theorem}
The nuclear Fr\'echet space ${\cal N}$ can be represented as 
$$
{\cal N=}\bigcap\limits_{p\in \N }{\cal H}_p\text{,} 
$$
where $\{{\cal H}_p$, $p\in {\N\}}$ is a family of Hilbert spaces such that
for all $p_1,p_2\in {\N}$ there exists $p\in {\N}$ such that the embeddings $%
{\cal H}_p\hookrightarrow {\cal H}_{p_1}$ and ${\cal H}_p\hookrightarrow 
{\cal H}_{p_2}$ are of Hilbert-Schmidt class. The topology of ${\cal N}$ is
given by the projective limit topology, i.e., the coarsest topology on ${\cal %
N}$ such that the canonical embeddings ${\cal N\hookrightarrow H}_p$ are
continuous for all $p\in \N $.
\end{theorem}

The Hilbertian norms on ${\cal H}_p$ are denoted by $\left| \cdot \right| _p$%
. Without loss of generality we always suppose that $\forall p\in \N %
,\forall \xi \in {\cal N}:\left| \xi \right| \leq \left| \xi \right| _p$ and
that the system of norms is ordered, i.e., $\left| \cdot \right| _p$ $\leq
\left| \cdot \right| _q$ if $p<q$. By general duality theory the dual space $%
{\cal N}^{\prime }\,$ can be written as 
$$
{\cal N}^{\prime }=\bigcup\limits_{p\in \N }{\cal H}_{-p}\text{.} 
$$
with inductive limit topology $\tau _{ind}$ by using the dual family of
spaces $\{{\cal H}_{-p}:={\cal H}_p^{\prime },\ p\in {\N\}}$. The inductive
limit topology (w.r.t. this family) is the finest topology on ${\cal N}%
^{\prime }$ such that the embeddings ${\cal H}_{-p}\hookrightarrow {\cal N}%
^{\prime }$ are continuous for all $p\in \N $. It is convenient to denote
the norm on ${\cal H}_{-p}$ by $\left| \cdot \right| _{-p}$. Let us mention
that in our setting the topology $\tau _{ind}$ coincides with the Mackey
topology $\tau ({\cal N}^{\prime },{\cal N})$ and the strong topology $\beta
({\cal N}^{\prime },{\cal N})$. Further note that the dual pair $\langle 
{\cal N}^{\prime },{\cal N}\rangle $ is reflexive if ${\cal N}^{\prime }$ is
equipped with $\beta ({\cal N}^{\prime },{\cal N})$. In addition we have
that convergence of sequences is equivalent in $\beta ({\cal N}^{\prime },%
{\cal N})$ and the weak topology $\sigma ({\cal N}^{\prime },{\cal N})$, see
e.g., \cite[Appendix 5]{HKPS93}.

Further we want to introduce the notion of tensor power of a nuclear space.
The simplest way to do this is to start from usual tensor powers ${\cal H}%
_p^{\otimes n}\ ,\ n\in \N $ of Hilbert spaces. Since there is no danger of
confusion we will preserve the notation $\left| \cdot \right| _p$ and $%
\left| \cdot \right| _{-p}$ for the norms on ${\cal H}_p^{\otimes n}$ and $%
{\cal H}_{-p}^{\otimes n}$ respectively. Using the definition%
$$
{\cal N}^{\otimes n}:=\ \stackunder{p\in \N }{\rm pr\ lim}{\cal H}_p^{\otimes
n} 
$$
one can prove \cite{Sch71} that ${\cal N}^{\otimes n}$ is a nuclear space
which is called the $n^{th}$ tensor power of ${\cal N.}$ The dual space of $%
{\cal N}^{\otimes n}$ can be written%
$$
\left( {\cal N}^{\otimes n}\right) ^{\prime }=\ \stackunder{p\in \N }{\rm ind\
lim}{\cal H}_{-p}^{\otimes n} 
$$

Most important for the applications we have in mind is the following 'kernel
theorem', see e.g., \cite{BeKo88}.

\begin{theorem}
\label{KernelTh}Let $\xi _1,...,\xi _n\mapsto F_n\left( \xi _1,...,\xi
_n\right) $ be an n-linear form on ${\cal N}^{\otimes n}$ which is ${\cal H}%
_p$-continuous , i.e.,%
$$
|F_n\left( \xi _1,...,\xi _n\right) |\leq C\prod_{k=1}^n|\xi _k|_p 
$$
for some $p\in {\N}$ and $C>0.$\\Then for all $p^{\prime }>p$ such that the
embedding $i_{p^{\prime },p}:{\cal H}_{p^{\prime }}\hookrightarrow {\cal H}_p
$ is Hilbert-Schmidt there exist a unique $\Phi ^{(n)}\in {\cal H}%
_{-p^{\prime }}^{\otimes n}$ such that 
$$
F_n\left( \xi _1,\ldots ,\xi _n\right) =\langle \Phi ^{(n)},\xi _1\otimes
\cdots \otimes \xi _n\rangle \ ,\quad \text{ }\xi _1,...,\xi _n\in {\cal N} 
$$
and the following norm estimate holds%
$$
\left| \Phi ^{(n)}\right| _{-p^{\prime }}\leq C\;\left\| i_{p^{\prime
},p}\right\| _{HS}^n 
$$
using the Hilbert-Schmidt norm of $i_{p^{\prime },p}$.
\end{theorem}

\begin{corollary}
\label{KernelCor}Let $\xi _1,...,\xi _n\mapsto F\left( \xi _1,...,\xi
_n\right) $ be an $n$-linear form on ${\cal N}^{\otimes n}$ which is ${\cal H%
}_{-p}$-continuous, i.e.,%
$$
|F_n\left( \xi _1,\ldots ,\xi _n\right) |\leq C\prod_{k=1}^n|\xi _k|_{-p} 
$$
for some $p\in {\N}$ and $C>0$.\\Then for all $p^{\prime }<p$ such that the
embedding $i_{p,p^{\prime }}:{\cal H}_p\hookrightarrow {\cal H}_{p^{\prime }}
$ is Hilbert-Schmidt there exist a unique $\Phi ^{(n)}\in {\cal H}%
_{p^{\prime }}^{\otimes n}$ such that 
$$
F_n\left( \xi _1,...,\xi _n\right) =\langle \Phi ^{(n)},\xi _1\otimes \cdots
\otimes \xi _n\rangle ,\quad \text{ }\xi _1,...,\xi _n\in {\cal N} 
$$
and the following norm estimate holds%
$$
\left| \Phi ^{(n)}\right| _{p^{\prime }}\leq C\;\left\| i_{p,p^{\prime
}}\right\| _{HS}^n\text{ .} 
$$
\end{corollary}

If in Theorem \ref{KernelTh} (and in Corollary \ref{KernelCor} respectively
) we start from a symmetric $n$-linear form $F_n$ on ${\cal N}^{\otimes n}$
i.e., $F_n(\xi _{\pi _1},\ldots ,\xi _{\pi _n})=F_n\left( \xi _1,\ldots ,\xi
_n\right) $ for any permutation $\pi $, then the corresponding kernel $\Phi
^{(n)}$ can be localized in ${\cal H}_{p^{\prime }}^{\hat \otimes n}\subset 
{\cal H}_{p^{\prime }}^{\otimes n}$ (the n$^{th}$ symmetric tensor power of
the Hilbert space ${\cal H}_{p^{\prime }}$). For $f_1,\ldots ,f_n\in {\cal H}
$ let $\hat \otimes $ also denote the symmetrization of the tensor product%
$$
f_1\hat \otimes \cdots \hat \otimes f_n:=\frac 1{n!}\sum_\pi f_{\pi
_1}\otimes \cdots \otimes f_{\pi _n}\ , 
$$
where the sum extends over all permutations of $n$ letters. All the above
quoted theorems also hold for complex spaces, in particular the complexified
space ${\cal N}_{\Ckl }$. By definition an element $\theta \in $ ${\cal N}_{%
\Ckl }$ decomposes into $\theta =\xi +i\eta \ ,\ \xi ,\eta \in {\cal N}$. If
we also introduce the corresponding complexified Hilbert spaces ${\cal H}_{p,%
\Ckl }$ the inner product becomes%
$$
(\theta _1,\theta _2)_{{\cal H}_{p,\Ckkl }}=(\theta _1,\bar \theta _2)_{%
{\cal H}_p}=(\xi _1,\xi _2)_{{\cal H}_p}+(\eta _1,\eta _2)_{{\cal H}%
_p}+i(\eta _1,\xi _2)_{{\cal H}_p}-i(\xi _1,\eta _2)_{{\cal H}_p} 
$$
for $\theta _1,\theta _2\in {\cal H}_{p,\Ckl },\ \theta _1=\xi _1+i\eta _1\
,\theta _2=\xi _2+i\eta _2\ ,\ \xi _1,\xi _2,\eta _1,\eta _2\in {\cal H}_p$.
Thus we have introduced a nuclear triple%
$$
{\cal N}_{\Ckl }^{\hat \otimes n}\subset {\cal H}_{\Ckl }^{\hat \otimes
n}\subset \left( {\cal N}_{\Ckl }^{\hat \otimes n}\right) ^{\prime } 
$$
We also want to introduce the (Boson or symmetric) Fock space $\Gamma ({\cal %
H})$ of ${\cal H}$ by%
$$
\Gamma ({\cal H})=\bigoplus_{n=0}^\infty {\cal H}_{\Ckl }^{\hat \otimes n} 
$$
with the convention ${\cal H}_{\Ckl }^{\hat \otimes 0}:=\C $ and the
Hilbertian norm%
$$
\left\| \vec \varphi \right\| _{\Gamma ({\cal H})}^2=\sum_{n=0}^\infty
n!\;\left| \varphi ^{(n)}\right| ^2\ ,\quad \vec \varphi =\left\{ \varphi
^{(n)}\;\Big|\;n\in \N _0\right\} \in \Gamma ({\cal H})\;. 
$$

%% file: holo.tex
\LaTeXparent{nga4.tex}

We shall collect some facts from the theory of holomorphic functions in
locally convex topological vector spaces ${\cal E}$ (over the complex field $%
{\C}$), see e.g., \cite{Di81}. Let ${\cal L}({\cal E}^n)$ be the space of
n-linear mappings from ${\cal E}^n$ into ${\C}$ and ${\cal L}_s({\cal E}^n)$
the subspace of symmetric n-linear forms. Also let ${\sl P}^n({\cal E})$
denote the n-homogeneous polynomials on ${\cal E}$. There is a linear
bijection ${\cal L}_s({\cal E}^n)\ni A\longleftrightarrow \widehat{A}\in 
{\sl P}^n({\cal E})$. Now let ${\cal U}\subset {\cal E}$ be open and
consider a function $G:{\cal U}\rightarrow {\C}$.

$G$ is said to be {\bf G-holomorphic} if for all $\theta _0\in {\cal U}$ and
for all $\theta \in {\cal E}$ the mapping from ${\C}$ to ${\C :}$ $\lambda
\rightarrow G(\theta _0+\lambda \theta )$ is holomorphic in some
neighborhood of zero in ${\C}$. If $G$ is G-holomorphic then there exists
for every $\eta \in {\cal U}$ a sequence of homogeneous polynomials $\frac
1{n!}\widehat{{\rm d}^n G(\eta )}$ such that 
$$
G(\theta +\eta )=\sum\limits_{n=0}^\infty \frac 1{n!}\widehat{{\rm d}%
^nG(\eta )}(\theta ) 
$$
for all $\theta $ from some open set ${\cal V}\subset {\cal U}$. $G$ is said
to be {\bf holomorphic}, if for all $\eta $ in ${\cal U}$ there exists an
open neighborhood ${\cal V}$ of zero such that $\sum\limits_{n=0}^\infty
\frac 1{n!}\widehat{{\rm d}^nG(\eta )}(\theta )$ converges uniformly on $%
{\cal V}$ to a continuous function. We say that $G$ is holomorphic at $%
\theta _0$ if there is an open set ${\cal U}$ containing $\theta _0$ such
that $G$ is holomorphic on ${\cal U}$. The following proposition can be
found e.g., in \cite{Di81}.

\begin{proposition}
\label{GHolLocB} $G$ is holomorphic if and only if it is G-holomorphic and
locally bounded.
\end{proposition}

\noindent Let us explicitly consider a function holomorphic at the point $%
0\in {\cal E}={\cal N}_{\Ckl }$, then

1) there exist $p$ and $\varepsilon >0$ such that for all $\xi _0\in {\cal N}%
_{\Ckl }$ with $\left| \xi _0\right| _p\leq \varepsilon $ and for all $\xi
\in {\cal N}_{\Ckl }$ the function of one complex variable $\lambda
\rightarrow G(\xi _0+\lambda \xi )$ is analytic at $0\in {\C}$, and

2) there exists $c>0$ such that for all $\xi \in {\cal N}_{\Ckl }$ with $%
\left| \xi \right| _p\leq \varepsilon $ : $\left| G(\xi )\right| \leq c$.

\noindent As we do not want to discern between different restrictions of one
function, we consider germs of holomorphic functions, i.e., we identify $F$
and $G$ if there exists an open neighborhood ${\cal U}:0\in {\cal U}\subset 
{\cal N}_{\Ckl }$ such that $F(\xi )=G(\xi )$ for all $\xi \in {\cal U}$.
Thus we define ${\rm Hol}_0({\cal N}_{\Ckl })$ as the algebra of germs of
functions holomorphic at zero equipped with the inductive topology given by
the following family of norms%
$$
{\rm n}_{p,l,\infty }(G)=\sup _{\left| \theta \right| _p\leq 2^{-l}}\left|
G(\theta )\right| ,\quad p,l\in \N . 
$$
\bigskip

\begin{sloppypar}
Let use now introduce spaces of entire functions which will be useful later.
Let ${\cal E}_{2^{-l}}^k({\cal H}_{-p,\Ckl })$ denote the set of all entire
functions on ${\cal H}_{-p,\Ckl }$ of growth $k\in [1,2]$ and type $2^{-l},\
p,l\in \Z $. This is a linear space with norm 
$$
{\rm n}_{p,l,k}(\varphi )=\sup _{z\in {\cal H}_{-p,\Ckkl }}\left| \varphi
(z)\right| \exp \left( -2^{-l}|z|_{-p}^k\right) ,\qquad \varphi \in {\cal E}%
_{2^{-l}}^k({\cal H}_{-p,\Ckl }) 
$$
The space of entire functions on ${\cal N}_{\Ckl }^{\prime }$ of growth $k$
and minimal type is naturally introduced by 
$$
{\cal E}_{\min }^k({\cal N}_{\Ckl }^{\prime }):=\ \stackunder{p,l\in \N }{\rm pr\ lim\,}{\cal E}_{2^{-l}}^k({\cal H}_{-p,\Ckl })\ , 
$$
see e.g., \cite{Ou91}. We will also need the space of entire functions on $%
{\cal N}_{\Ckl }$ of growth $k$ and finite type:%
$$
{\cal E}_{\max }^k({\cal N}_{\Ckl }):=\ \stackunder{p,l\in \N }{\rm ind\ lim}%
{\cal E}_{2^l}^k({\cal H}_{p,\Ckl })\ . 
$$
In the following we will give an equivalent description of ${\cal E}_{\min
}^k({\cal N}_{\Ckl }^{\prime })$ and ${\cal E}_{\max }^k({\cal N}_{\Ckl })$.
Cauchy's inequality and Corollary \ref{KernelCor} allow to write the Taylor
coefficients in a convenient form. Let $\varphi \in {\cal E}_{\min }^k({\cal %
N}_{\Ckl }^{\prime })$ and $z\in {\cal N}_{\Ckl }^{\prime }$, then there
exist kernels $\varphi ^{(n)}\in {\cal N}_{\Ckl }^{\hat \otimes n}$ such
that 
$$
\langle z^{\otimes n},\varphi ^{(n)}\rangle =\frac 1{n!}\widehat{{\rm d}^n\varphi (0)}(z) 
$$
i.e., 
\begin{equation}
\label{phi(z)}\varphi (z)=\sum_{n=0}^\infty \langle z^{\otimes n},\varphi
^{(n)}\ \rangle . 
\end{equation}
This representation allows to introduce a nuclear topology on ${\cal E}%
_{\min }^k({\cal N}_{\Ckl }^{\prime })$, see \cite{Ou91} for details. Let 
{\rm E}$_{p,q}^\beta $ denote the space of all functions of the form (\ref
{phi(z)}) such that the following Hilbertian norm 
\begin{equation}
\label{3StrichNorm}\lnorm  \varphi \rnorm  _{p,q,\beta
}^2:=\sum_{n=0}^\infty (n!)^{1+\beta }2^{nq}\left| \varphi ^{(n)}\right|
_p^2\;,\quad p,q\in \N  
\end{equation}
is finite for $\beta \in [0,1]$. (By $\left| \varphi ^{(0)}\right| _p$ we 
simply mean the complex modulus for
all $p$.)  The space {\rm E}$_{-p-,q}^{-\beta }$ with
the norm $\lnorm  \varphi \rnorm  _{-p,-q,-\beta }$ is defined analogously.
\end{sloppypar}

\begin{theorem}
\label{Ekminprlim}The following topological identity holds: 
$$
\stackunder{p,q\in \N }{\rm pr\ lim}\;{\rm E}_{p,q}^\beta ={\cal E}_{\min
}^{\frac 2{1+\beta }}({\cal N}_{\Ckl }^{\prime })\quad . 
$$
\end{theorem}

The proof is an immediate consequence of the following two lemmata which
show that the two systems of norms are in fact equivalent.

\begin{lemma}
\label{nplk3Strich}Let $\varphi \in ${\rm E}$_{p,q}^\beta $ then $\varphi
\in {\cal E}_{2^{-l}}^{\frac 2{1+\beta }}({\cal H}_{-p,\Ckl })$ for $l=\frac
q{1+\beta }$. Moreover 
\begin{equation}
\label{nplk3StrichNorm}{\rm n}_{p,l,k}(\varphi )\leq \lnorm  \varphi \rnorm  %
_{p,q,\beta }\ ,\ \ k=\tfrac 2{1+\beta }\ .
\end{equation}
\end{lemma}

\TeXButton{Proof}{\proof} We look at the convergence of the series $\varphi
(z)=\sum_{n=0}^\infty \langle z^{\otimes n},\varphi ^{(n)}\ \rangle \ $, $%
z\in {\cal H}_{-p,\Ckl }\ ,\ \varphi ^{(n)}\in {\cal H}_{p,\Ckl }$ if $%
\sum_{n=0}^\infty (n!)^{1+\beta }2^{nq}|\varphi ^{(n)}|_p^2=\lnorm 
\varphi \rnorm  _{p,q,\beta }^2$ is finite. The following estimate holds:%
\begin{eqnarray*}
\sum_{n=0}^\infty |\langle z^{\otimes n},\varphi ^{(n)}\ \rangle | %
& \leq & \left( \sum_{n=0}^\infty (n!)^{1+\beta }2^{nq}|\varphi %
^{(n)}|_p^2\right) ^{1/2}\left( \sum_{n=0}^\infty %
\frac 1{(n!)^{1+\beta }}2^{-nq}|z|_{-p}^{2n}\right) ^{1/2} %
\\& \leq & \lnorm  \varphi \rnorm  _{p,q,\beta }\cdot \left( %
\sum_{n=0}^\infty \left\{ \frac 1{n!}2^{-\frac{nq}{1+\beta }}%
|z|_{-p}^{\frac{2n}{1+\beta }}\right\} ^{1+\beta }\right) ^{1/2} %
\\& \leq & \lnorm  \varphi \rnorm  _{p,q,\beta }\left( %
\sum_{n=0}^\infty \frac 1{n!}2^{-\frac{nq}{1+\beta }}|z|_{-p}^{\frac{2n}{1+\beta }}\right) ^{(1+\beta )/2} \\& \leq & \lnorm  \varphi \rnorm %
_{p,q,\beta }\exp \left( %
2^{-\frac q{1+\beta }}|z|_{-p}^{\frac 2{1+\beta }}\right) . 
\end{eqnarray*}
\TeXButton{End Proof}{\endproof}

\begin{lemma}
\label{3Strichnplk}For any $p^{\prime },q\in \N $ there exist $p,l\in \N $
such that 
$$
{\cal E}_{2^{-l}}^{\frac 2{1+\beta }}({\cal H}_{-p,\Ckl })\subset {\rm E}%
_{p^{\prime },q}^\beta  
$$
i.e., there exists a constant $C>0$ such that 
$$
\lnorm  \varphi \rnorm  _{p^{\prime },q,\beta }\leq C\;{\rm n}%
_{p,l,k}(\varphi ),\quad \varphi \in {\cal E}_{2^{-l}}^k({\cal H}_{-p,\Ckl %
}),\quad k=\tfrac 2{1+\beta }. 
$$
\end{lemma}

\TeXButton{Remark}{\remark } More precisely we will prove the following:\ If $%
\varphi \in {\cal E}_{2^{-l}}^k({\cal H}_{-p,\Ckl })$ then $\varphi \in $%
{\rm E}$_{p^{\prime },q}^\beta $ for $k=\frac 2{1+\beta }$ and $\rho
:=2^{q-2l/k}k^{2/k}e^2\left\| i_{p^{\prime },p}\right\| _{HS}^2<1$ (in
particular this requires $p^{\prime }>p$ to be such that the embedding $%
i_{p^{\prime },p}:{\cal H}_{p^{\prime }}\hookrightarrow {\cal H}_p$ is
Hilbert-Schmidt). \\Moreover the following bound holds 
\begin{equation}
\label{3Normnplk}\lnorm  \varphi \rnorm  _{p^{\prime },q,\beta }\leq {\rm n}%
_{p,l,k}(\varphi )\cdot \left( 1-\rho \right) ^{-1/2}\ . 
\end{equation}

\TeXButton{Proof}{\proof} The assumption $\varphi \in {\cal E}_{2^{-l}}^k(%
{\cal H}_{-p,\Ckl })$ implies a bound of the growth of $\varphi :$%
$$
|\varphi (z)|\leq {\rm n}_{p,l,k}(\varphi )\exp (2^{-l}|z|_{-p}^k)\ . 
$$
For each $\rho >0\ ,\ z\in {\cal H}_{-p,\Ckl }$ the Cauchy inequality from
complex analysis \cite{Di81} gives%
$$
\left| \frac 1{n!}\widehat{{\rm d}^n\varphi (0)}(z)\right| \leq {\rm n}%
_{p,l,k}(\varphi )\rho ^{-n}\exp (\rho ^k2^{-l})\;|z|_{-p}^n\ . 
$$
By polarization \cite{Di81} it follows for $z_1,\ldots ,z_n\in {\cal H}_{-p,\Ckl }$%
$$
\left| \frac 1{n!}{\rm d}^n\varphi (0)(z_1,\ldots ,z_n)\right| \leq {\rm n}%
_{p,l,k}(\varphi )\frac 1{n!}\left( \frac n\rho \right) ^n\exp (\rho
^k2^{-l})\prod_{k=1}^n|z_k|_{-p}\ . 
$$
For $p^{\prime }>p$ such that $\left\| i_{p^{\prime },p}\right\| _{HS}$ is
finite, an application of the kernel theorem guarantees the existence of
kernels $\varphi ^{(n)}\in {\cal H}_{p^{\prime },\Ckl }^{\hat \otimes n}$
such that 
$$
\varphi (z)=\sum_{n=0}^\infty \langle z^{\hat \otimes n},\varphi ^{(n)}\
\rangle 
$$
with the bound%
$$
\left| \varphi ^{(n)}\right| _{p^{\prime }}\leq {\rm n}_{p,l,k}(\varphi
)\frac 1{n!}\left( \frac n\rho \left\| i_{p^{\prime },p}\right\|
_{HS}\right) ^n\exp (\rho ^k\cdot 2^{-l})\ . 
$$
We can optimize the bound with the choice of an $n$-dependent $\rho $.
Setting $\rho ^k=2^ln/k$ we obtain 
\begin{eqnarray*}
\left| \varphi ^{(n)}\right| _{p^{\prime }} & \leq & {\rm n}_{p,l,k}(\varphi %
)\frac 1{n!}n^{n(1-1/k)}\left( \tfrac 1k2^l\right) ^{-n/k}\left\| %
i_{p^{\prime },p}\right\| _{HS}^ne^{n/k} %
\\& \leq & {\rm n}_{p,l,k}(\varphi )\;(n!)^{-1/k}\left\{ (k2^{-l})^{1/k}e %
\left\| i_{p^{\prime },p}\right\| _{HS}\right\} ^n\ , 
\end{eqnarray*}
where we used $n^n\leq n!\,e^n$ in the last estimate. Now choose $\beta \in
[0,1]$ such that $k=\frac 2{1+\beta }$ to estimate the following norm:%
\begin{eqnarray*}
\lnorm  \varphi \rnorm  _{p^{\prime },q,\beta }^2 %
& \leq & {\rm n}%
_{p,l,k}^2(\varphi )\sum_{n=0}^\infty (n!)^{1+\beta -\frac 2k}2^{qn}\left\{ %
(k2^{-l})^{1/k}e\left\| i_{p^{\prime },p}\right\| _{HS}\right\} ^{2n} %
\\& \leq & {\rm n}_{p,l,k}^2(\varphi )\left( 1-2^q\left\{ (k2^{-l})^{1/k}e %
\left\| i_{p^{\prime },p}\right\| _{HS}\right\} ^2\right) ^{-1} 
\end{eqnarray*}
for sufficiently large $l$. This completes the proof.\TeXButton{End Proof}
{\endproof}\bigskip\ 

Analogous estimates for these systems of norms also hold if $\beta ,p,q,l$
become negative. This implies the following theorem. For related results see
e.g., \cite[Prop.8.6]{Ou91}.

\begin{theorem}
\label{indlimEkmax} \hfill \\If $\beta \in [0,1)$ then the following
topological identity holds:%
$$
\stackunder{p,q\in \N }{\rm ind\ lim}\ {\rm E}_{-p,-q}^{-\beta }={\cal E}%
_{\max }^{2/(1-\beta )}({\cal N}_{\Ckl }). 
$$
If $\beta =1$ we have%
$$
\stackunder{p,q\in \N }{\rm ind\ lim}\ {\rm E}_{-p,-q}^{-1}={\rm Hol}_0(%
{\cal N}_{\Ckl })\ . 
$$
\end{theorem}

\noindent This theorem and its proof will appear in the context of section 
\ref{Characterization}. The characterization of distributions in infinite
dimensional analysis is strongly related to this theorem. From this point of
view it is natural to postpone its proof to section \ref{Characterization}.

%% file: measures.tex
\LaTeXparent{nga4.tex}

\section{Measures on linear topological spaces}

To introduce probability measures on the vector space ${\cal N}^{\prime }$,
we consider ${\cal C}_\sigma ({\cal N^{\prime }})$ the $\sigma $-algebra
generated by cylinder sets on ${\cal N}^{\prime }$, which coincides with the
Borel $\sigma $-algebras ${\cal B}_\sigma ({\cal N}^{\prime })$ and ${\cal B}%
_\beta ({\cal N}^{\prime })$ generated by the weak and strong topology on $%
{\cal N}^{\prime }$ respectively. Thus we will consider this $\sigma $%
-algebra as the {\em natural} $\sigma $-algebra on ${\cal N}^{\prime }$.
Detailed definitions of the above notions and proofs of the mentioned
relations can be found in e.g., \cite{BeKo88}.

We will restrict our investigations to a special class of measures $\mu $ on 
${\cal C}_\sigma ({\cal N^{\prime }})$, which satisfy two additional
assumptions. The first one concerns some analyticity of the Laplace
transformation%
$$
l_\mu (\theta )=\int_{{\cal N^{\prime }}}\exp \left\langle x,\theta
\right\rangle \ {\rm d}\mu (x)=:\E _\mu (\exp \left\langle \cdot ,\theta
\right\rangle )\text{ , }\theta \in {\cal N}_{\Ckl }\text{.} 
$$
Here we also have introduced the convenient notion of expectation $\E _\mu $
of a $\mu $-integrable function.

\TeXButton{Assumption}{\assumption} 1 \quad The measure $\mu $ has an analytic
Laplace transform in a neighborhood of zero. That means there exists an open
neighborhood ${\cal U}\subset {\cal N}_{\Ckl }$ of zero, such that $l_\mu $
is holomorphic on ${\cal U}$, i.e., $l_\mu \in {\rm Hol}_0({\cal N}_{\Ckl })$
. This class of {\em analytic measures} is denoted by ${\cal M}_a({\cal N}%
^{\prime }).$\bigskip\ 

\noindent An equivalent description of analytic measures is given by the
following lemma.

\begin{lemma}
\label{equiLemma}The following statements are equivalent\medskip\ \\{\bf 1)}%
\quad $\mu \in {\cal M}_a({\cal N}^{\prime })$\smallskip\\{\bf 2)}$\quad \D%
\exists p_\mu \in {\N},\quad \exists C>0:\qquad \left| \int_{{\cal N}%
^{\prime }}\ \langle x,\theta \rangle ^n\;{\rm d}\mu (x)\right| \leq
n!\,C^n\left| \theta \right| _{p_\mu }^n\;,\quad \theta \in {\cal H}_{p_\mu ,%
\Ckl }$\smallskip\\{\bf 3)}$\quad \D\exists p_\mu ^{\prime }\in {\N},\quad
\exists \varepsilon _\mu >0:\ \qquad \int_{{\cal N}^{\prime }}\ \exp
(\varepsilon _\mu \left| x\right| _{-p_\mu ^{\prime }})\,{\rm d}\mu
(x)<\infty $
\end{lemma}

\noindent \TeXButton{Proof}{\proof}The proof can be found in \cite{KoSW94}.
We give its outline in the following. The only non-trivial step is the proof
of 2)$\Rightarrow $3).

\noindent By polarization \cite{Di81} 2) implies 
\begin{equation}
\label{PolBound}\left| \int_{{\cal N}^{\prime }}\langle x^{\otimes
n},\tbigotimes_{j=1}^n\xi _j\rangle \ {\rm d}\mu (x)\right| \leq
n!\;C^n\prod_{j=1}^n\left| \xi _j\right| _{p_\mu }\ ,\quad \xi _j\in {\cal H}%
_{p^{\prime }} 
\end{equation}
for a (new) constant $C>0$. Choose $p^{\prime }>p_\mu $ such that the
embedding $i_{p^{\prime },p_\mu }:{\cal H}_{p^{\prime }}\rightarrow {\cal H}%
_{p_\mu }$ is of Hilbert-Schmidt type. Let $\left\{ e_k,\ k\in {\N}\right\}
\subset {\cal N}$ be an orthonormal basis in ${\cal H}_{p^{\prime }}$. Then $%
\left| x\right| _{-p^{\prime }}^2=\sum\limits_{k=1}^\infty \left\langle
x,e_k\right\rangle ^2$, $x\in {\cal H}_{-p^{\prime }}$. We will first
estimate the moments of even order%
$$
\int_{{\cal N}^{\prime }}\left| x\right| _{-p^{\prime }}^{2n}\ {\rm d}\mu
(x)=\sum\limits_{k_1=1}^\infty \cdots \sum\limits_{k_n=1}^\infty \int_{{\cal %
N}^{\prime }}\left\langle x,e_{k_1}\right\rangle ^2\ \cdots \left\langle
x,e_{k_n}\right\rangle ^2\ {\rm d}\mu (x)\ \text{,} 
$$
where we changed the order of summation and integration by a monotone
convergence argument. Using the bound (\ref{PolBound}) we have%
\begin{eqnarray*}
\int_{{\cal N}^{\prime }}\left| x\right| _{-p^{\prime }}^{2n}\ {\rm d}\mu(x)%
&\leq & \ C^{2n}\ (2n)!\sum\limits_{k_1=1}^\infty \cdots%
\sum\limits_{k_n=1}^\infty \left| e_{k_1}\right| _{p_\mu}^2\cdots \left|%
e_{k_n}\right| _{p_\mu}^2 %
\\&=& \ C^{2n}\ (2n)!\left( \sum\limits_{k=1}^\infty \left| e_k\right|%
_{p_\mu}^2\right) ^n %
\\&=& \ \left( C\cdot \left\| i_{p^{\prime },p_\mu}\right\| _{HS}%
\right) ^{2n}(2n)!%
\end{eqnarray*}
because%
$$
\sum\limits_{k=1}^\infty \left| e_k\right| _{p_\mu}^2=\left\| i_{p^{\prime
},p_\mu }\right\| _{HS}^2\text{ .} 
$$
The moments of arbitrary order can now be estimated by the Schwarz inequality%
\begin{eqnarray*}
\int \left| x\right| _{-p^{\prime }}^n\ {\rm d}\mu (x) %
& \leq & \sqrt{\mu ({\cal N}^{\prime })}\left( \int \left| x\right| _{-p}^{2n}%
\ {\rm d}\mu (x)\right)^{\frac 12} %
\\&\leq & \sqrt{\mu ({\cal N}^{\prime })}\left( C\left\| %
i_{p^{\prime },p_\mu}\right\| _{HS}\right) ^n\sqrt{(2n)!} %
\\&\leq & \sqrt{\mu ({\cal N}^{\prime })}\left( 2 C\left\| %
i_{p^{\prime },p_\mu}\right\| _{HS}\right) ^nn! 
\end{eqnarray*}
since $(2n)!\leq 4^n(n!)^2$ . \\Choose $\varepsilon <$ $\left( \ 2C\left\|
i_{p^{\prime },p_\mu }\right\| _{HS}\right) ^{-1}$ then%
\begin{eqnarray} \label{exponent}
\int e^{\varepsilon \left| x\right| _{-p^{\prime }}}{\rm d}\mu%
(x) &=& \sum_{n=0}^\infty \frac{\varepsilon ^n}{n!}\int \left| x\right|%
_{-p^{\prime }}^n\ {\rm d}\mu (x) \nonumber %
\\&\leq & \sqrt{\mu ({\cal N}^{\prime })}%
\ \sum_{n=0}^\infty \left( \varepsilon \ 2 C%
\left\| i_{p^{\prime },p_\mu }\right\| _{HS}\right) ^n<\infty 
\end{eqnarray}
Hence the lemma is proven.\TeXButton{End Proof}{\endproof}\bigskip\ \ 

For $\mu \in {\cal M}_a({\cal N}^{\prime })$ the estimate in statement 2 of
the above lemma allows to define the moment kernels ${\rm M}_n^\mu \in (%
{\cal N}^{\hat \otimes n})^{\prime }.$ This is done by extending the above
estimate by a simple polarization argument and applying the kernel theorem.
The kernels are determined by%
$$
l_\mu (\theta )=\sum_{n=0}^\infty \frac 1{n!}\langle {\rm M}_n^\mu ,\theta
^{\otimes n}\rangle 
$$
or equivalently%
$$
\langle {\rm M}_n^\mu ,\theta _1\hat \otimes \cdots \hat \otimes \theta
_n\rangle =\left. \frac{\partial ^n}{\partial t_1\cdots \partial t_n}l_\mu
(t_1\theta _1+\cdots +t_n\theta _n)\right| _{t_1=\cdots =t_n=0}\ . 
$$
Moreover, if $p>p_\mu $ is such that embedding $i_{p,p_\mu }:{\cal H}%
_p\hookrightarrow {\cal H}_{p_\mu }$ is Hilbert-Schmidt then 
\begin{equation}
\label{MnmuNorm}\left| {\rm M}_n^\mu \right| _{-p}\leq \left( nC\left\|
i_{p,p_\mu }\right\| _{HS}\right) ^n\leq {n!}\ \left( eC\left\| i_{p,p_\mu
}\right\| _{HS}\right) ^n\ . 
\end{equation}

\begin{definition}
\label{P(N)}A function $\varphi :{\cal N}^{\prime }\rightarrow {\C}$ of the
form $\varphi (x)=\sum_{n=0}^N\langle x^{\otimes n},\varphi ^{(n)}\rangle $, 
$x\in {\cal N}^{\prime }$, $N\in {\N,}$ is called a continuous polynomial
(short $\varphi \in {\cal P}({\cal N}^{\prime })$ ) iff $\varphi ^{(n)}\in 
{\cal N}_{\Ckl }^{\hat \otimes n}$, $\forall n\in {\N}_0=\N \cup \{0\}$.
\end{definition}

Now we are ready to formulate the second assumption: \bigskip\ 

\TeXButton{Assumption}{\assumption} 2 \quad For all $\varphi \in {\cal P}({\cal N}%
^{\prime })$ with $\varphi =0$ $\mu $-almost everywhere we have $\varphi
\equiv 0$. In the following a measure with this property will be called {\em %
non-degenerate}. \bigskip\ 

\TeXButton{Note}{\note } Assumption 2 is equivalent to:\\Let $\varphi \in {\cal P}({\cal N}%
^{\prime })$ with $\int_A\varphi \,{\rm d}\mu =0$ for all $A\in {\cal C}%
_\sigma ({\cal N}^{\prime })$ then $\varphi \equiv 0$.\\A sufficient
condition can be obtained by regarding admissible shifts of the measure $\mu 
$. If $\mu (\cdot +\xi )$ is absolutely continuous with respect to $\mu $
for all $\xi \in {\cal N,}$ i.e., there exists the Radon-Nikodym derivative$%
{\cal \quad }$%
$$
\rho _\mu (\xi ,x)=\frac{{\rm d}\mu (x+\xi )}{{\rm d}\mu (x)}\ ,\quad x\in 
{\cal N}^{\prime }{\rm \;,} 
$$
Then we say that $\mu $ is ${\cal N}${\em --quasi-invariant} see e.g., \cite
{GV68,Sk74}. This is sufficient to ensure Assumption 2, see e.g., \cite
{KoTs91,BeKo88}.\bigskip\ 

\example In Gaussian Analysis (especially White Noise Analysis) the Gaussian
measure $\gamma _{{\cal H}}$ corresponding to the Hilbert space ${\cal H}$
is considered. Its Laplace transform is given by 
$$
l_{\gamma _{{\cal H}}}(\theta )=e^{\frac 12\langle \theta ,\theta \rangle }\
,\qquad \theta \in {\cal N}_{\Ckl }\ , 
$$
hence $\gamma _{{\cal H}}\in {\cal M}_a({\cal N}^{\prime })$. It is well
known that $\gamma _{{\cal H}}$ is ${\cal N}$--quasi-invariant (moreover $%
{\cal H}$--quasi-invariant) see e.g., \cite{Sk74,BeKo88}. Due to the
previous note $\gamma _{{\cal H}}$ satisfies also Assumption 2.\bigskip\ 

\example  {\it (Poisson measures)} \smallskip \\ Let use consider the
classical (real) Schwartz triple 
$$
{\cal S}(\R )\subset L^2(\R )\subset {\cal S}^{\prime }(\R )\,. 
$$
The Poisson white noise measure $\mu _p$ is defined as a probability measure
on ${\cal C}_\sigma ({\cal S}^{\prime }(\R ))$ with the Laplace transform%
$$
l_{\mu _p}(\theta )=\exp \left\{ \int_{\R }(e^{\theta (t)}-1)\;{\rm d}%
t\right\} =\exp \left\{ \langle e^\theta -1,1\rangle \right\} ,\qquad \theta
\in {\cal S}_{\Ckl }(\R )\,, 
$$
see e.g., \cite{GV68}. It is not hard to see that $l_{\mu _p}$ is a
holomorphic function on ${\cal S}_{\Ckl }(\R )$, so Assumption 1 is
satisfied. But to check Assumption 2, we need additional considerations.

First of all we remark that for any $\xi \in {\cal S}(\R )\,,\ \xi \neq 0$
the measures $\mu _p$ and $\mu _p(\cdot +\xi )$ are orthogonal (see \cite
{VGG75} for a detailed analysis). It means that $\mu _p$ is not ${\cal S}(\R %
)$-quasi-invariant and the note after Assumption 2 is not applicable now.

Let some $\varphi \in {\cal P}({\cal S}^{\prime }(\R ))\,,\,\varphi =0~~\mu
_p$-a.s. be given. We need to show that then $\varphi \equiv 0$. To this end
we will introduce a system of orthogonal polynomials in the space $L^2(\mu
_p)$ which can be constructed in the following way. The mapping%
$$
\theta (\cdot )\mapsto \alpha (\theta )(\cdot )=\log (1+\theta (\cdot ))\in 
{\cal S}_{\Ckl }(\R )\,,\quad \theta \in {\cal S}_{\Ckl }(\R ) 
$$
is holomorphic on a neighborhood ${\cal U}\subset {\cal S}_{\Ckl }(\R %
)\,,\,0\in {\cal U}$. Then%
$$
e_{\mu _p}^\alpha (\theta ;x)=\frac{e^{\langle \alpha (\theta ),x\rangle }}{%
l_{\mu _p}(\alpha (\theta ))}=\exp \{\langle \alpha (\theta ),x\rangle
-\langle \theta ,1\rangle \}\,,\quad \theta \in {\cal U}\,,\ x\in {\cal S}%
^{\prime }(\R ) 
$$
is a holomorphic function on ${\cal U}$ for any $\,x\in {\cal S}^{\prime }(%
\R )$. The Taylor decomposition and the kernel theorem (just as in
subsection \ref{AppellSec} below) give 
$$
e_{\mu _p}^\alpha (\theta ;x)=\sum_{n=0}^\infty \frac 1{n!}\langle \theta
^{\otimes n},C_n(x)\rangle \,, 
$$
where $C_n:{\cal S}^{\prime }(\R )\rightarrow {\cal S}^{\prime }(\R )^{\hat
\otimes n}$ are polynomial mappings. For $\varphi ^{(n)}\in {\cal S}_{\Ckl }(%
\R )^{\hat \otimes n}\,,\,n\in \N _0$, we define Charlier polynomials%
$$
x\mapsto C_n(\varphi ^{(n)};x)=\langle \varphi ^{(n)},C_n(x)\rangle \in \C %
\,,\ \,x\in {\cal S}^{\prime }(\R )\,. 
$$
Due to \cite{Ito88,IK88} we have the following orthogonality property:%
$$
\forall \varphi ^{(n)}\in {\cal S}_{\Ckl }(\R )^{\hat \otimes n}\,,\,\forall
\psi ^{(m)}\in {\cal S}_{\Ckl }(\R )^{\hat \otimes n} 
$$
$$
\int C_n(\varphi ^{(n)})C_m(\psi ^{(m)})\;{\rm d}\mu _p=\delta _{nm}
n!\langle \varphi ^{(n)},\psi ^{(n)}\rangle \,. 
$$
Now the rest is simple. Any continuous polynomial $\varphi $ has a uniquely
defined decomposition%
$$
\varphi (x)=\sum_{n=0}^N\langle \varphi ^{(n)},C_n(x)\rangle \,\,,\quad
\,x\in {\cal S}^{\prime }(\R )\,, 
$$
where $\varphi ^{(n)}\in {\cal S}_{\Ckl }(\R )^{\hat \otimes n}$. If $%
\varphi =0~\mu _p$-a.e. then 
$$
\left\| \varphi \right\| _{L^2(\mu _p)}^2=\sum_{n=0}^Nn!\,\langle \varphi
^{(n)},\overline{\varphi ^{(n)}}\rangle =0. 
$$
Hence $\varphi ^{(n)}=0\,,\,n=0\,,\,\ldots ,\,N$, i.e., $\varphi \equiv 0$.
So Assumption 2 is satisfied.

\section{Concept of distributions in infinite dimensional analysis\label
{Concept}}

In this section we will introduce a preliminary distribution theory in
infinite dimensional non-Gaussian analysis. We want to point out in advance
that the distribution space constructed here is in some sense too big for
practical purposes. In this sense section \ref{Concept} may be viewed as a
stepping stone to introduce the more useful structures in \S \ref
{testfunctions} and \S \ref{Distributions}.

We will choose ${\cal P}({\cal N}^{\prime })$ as our (minimal) test function
space. (The idea to use spaces of this type as appropriate spaces of test
functions is rather old see \cite{KMP65}. They also discussed in which sense
this space is ``minimal''.) First we have to ensure that ${\cal P}({\cal N}%
^{\prime })$ is densely embedded in $L^2(\mu )$. This is fulfilled because
of our assumption 1 \cite[Sec.§10 Th.1]{Sk74}. The space ${\cal P}({\cal N}%
^{\prime })$ may be equipped with various different topologies, but there
exists a natural one such that ${\cal P}({\cal N}^{\prime })$ becomes a
nuclear space \cite{BeKo88}. The topology on ${\cal P}({\cal N}^{\prime })$
is chosen such that is becomes isomorphic to the topological direct sum of
tensor powers ${\cal N}_{\Ckl }^{\hat \otimes n}$ see e.g., 
\cite[Ch II 6.1, Ch III 7.4]{Sch71}%
$$
{\cal P}({\cal N}^{\prime })\simeq \bigoplus_{n=0}^\infty {\cal N}_{\Ckl %
}^{\hat \otimes n}\text{ .} 
$$
via 
$$
\varphi (x)=\sum_{n=0}^\infty \left\langle x^{\otimes n},\varphi
^{(n)}\right\rangle \longleftrightarrow \vec \varphi =\left\{ \varphi
^{(n)}\;\Big|\;n\in \N _0\right\} . 
$$
Note that only a finite number of $\varphi ^{(n)}$ is non-zero. We will not
reproduce the full construction here, but we will describe the notion of
convergence of sequences this topology on ${\cal P}({\cal N}^{\prime })$.
For $\varphi \in {\cal P}({\cal N}^{\prime })$, $\varphi
(x)=\sum_{n=0}^{N(\varphi )}\left\langle x^{\otimes n},\varphi
^{(n)}\right\rangle $ let $p_n:{\cal P}({\cal N}^{\prime })\rightarrow {\cal %
N}_{\Ckl }^{\hat \otimes n}$ denote the mapping $p_n$ is defined by $%
p_n\varphi :=\varphi ^{(n)}.$ A sequence $\left\{ \varphi _j,\ j\in {\N}%
\right\} $ of smooth polynomials converges to $\varphi \in {\cal P}({\cal N}%
^{\prime })$ iff the $N(\varphi _j)\ $are bounded and $p_n\varphi _j%
\stackunder{n\rightarrow \infty }{\longrightarrow }p_n\varphi $ in ${\cal N}%
_{\Ckl }^{\hat \otimes n}$ for all $n\in {\N}$.

Now we can introduce the dual space ${\cal P}_\mu ^{\prime }({\cal N}%
^{\prime })$ of ${\cal P}({\cal N}^{\prime })$ with respect to $L^2(\mu )$.
As a result we have constructed the triple 
$$
{\cal P}({\cal N}^{\prime })\subset L^2(\mu )\subset {\cal P}_\mu ^{\prime }(%
{\cal N}^{\prime }) 
$$
The (bilinear) dual pairing $\langle \!\langle \cdot ,\cdot \rangle
\!\rangle _\mu $ between ${\cal P}_\mu ^{\prime }({\cal N}^{\prime })$ and $%
{\cal P}({\cal N}^{\prime })$ is connected to the (sesqui\-linear) inner
product on $L^2(\mu )$ by%
$$
\langle \!\langle \varphi ,\;\psi \rangle \!\rangle _\mu =(\varphi ,\;%
\overline{\psi })_{L^2(\mu )}\ ,\quad \varphi \in L^2(\mu ),\ \psi \in {\cal %
P}({\cal N}^{\prime })\text{ .} 
$$
Since the constant function 1 is in ${\cal P}({\cal N}^{\prime })$ we may
extend the concept of expectation from random variables to distributions;
for $\Phi \in {\cal P}_\mu ^{\prime }({\cal N}^{\prime })$ 
$$
\E _\mu (\Phi ):=\left\langle \!\left\langle \Phi ,1\right\rangle
\!\right\rangle _\mu \;. 
$$
The main goal of this section is to provide a description of ${\cal P}_\mu
^{\prime }({\cal N}^{\prime })$ , see Theorem \ref{PStrichRep} below. The
simplest approach to this problem seems to be the use of so called $\mu $%
-Appell polynomials.

\subsection{Appell polynomials associated to the measure
\texorpdfstring{$\mu $ }{ \textmu }
\label{AppellSec}}

Because of the holomorphy of $l_\mu $ and $l_\mu (0)=1$ there exists a
neighborhood of zero%
$$
{\cal U}_0=\left\{ \theta \in {\cal N}_{\Ckl }\ \Big|\ 2^{q_0}\left| \theta
\right| _{p_0}<1\right\} 
$$
$p_0,q_0\in {\N,}$ $p_0\geq p_\mu ^{\prime }$ , $2^{-q_0}\leq \varepsilon
_\mu $ ($p_\mu ^{\prime },\varepsilon _\mu $ from Lemma \ref{equiLemma})
such that $l_\mu (\theta )\neq 0$ for $\theta \in {\cal U}_0$ and the
normalized exponential 
\begin{equation}
\label{emy}e_\mu (\theta ;z)=\frac{e^{\left\langle z,\theta \right\rangle }}{%
l_\mu (\theta )}\text{ \quad for }\theta \in {\cal U}_0,\quad z\in {\cal N}_{%
\Ckl }^{\prime }\;, 
\end{equation}
is well defined. We use the holomorphy of $\theta \mapsto e_\mu (\theta ;z)$
to expand it in a power series in $\theta $ similar to the case
corresponding to the construction of one dimensional Appell polynomials \cite
{Bo76}. We have in analogy to \cite{AKS93,ADKS94}%
$$
e_\mu (\theta ;z)=\sum_{n=0}^\infty \frac 1{n!}\widehat{{\rm d}^ne_\mu (0,z)}%
(\theta ) 
$$
where $\widehat{{\rm d}^ne_\mu (0;z)}$ is an n-homogeneous continuous
polynomial. But since $e_\mu (\theta ;z)$ is not only G-holomorphic but
holomorphic we know that $\theta \rightarrow $ $e_\mu (\theta ;z)$ is also
locally bounded. Thus Cauchy's inequality for Taylor series \cite{Di81} may
be applied, $\rho \leq 2^{-q_0}$ , $p\geq p_0$%
\begin{equation}
\label{Cauchyemy}\left| \frac 1{n!}\widehat{{\rm d}^ne_\mu (0;z)}(\theta
)\right| \leq \frac 1{\rho ^n}\sup \limits_{\left| \theta \right| _p=\rho
}\left| e_\mu (\theta ;z)\right| \left| \theta \right| _p^n\leq \frac 1{\rho
^n}\sup \limits_{\left| \theta \right| _p=\rho }\frac 1{l_\mu (\theta
)}e^{\rho \left| z\right| _{-p}}\left| \theta \right| _p^n 
\end{equation}
if $z\in {\cal H}_{-p,\Ckl }$. This inequality extends by polarization \cite
{Di81} to an estimate sufficient for the kernel theorem. Thus we have a
representation $\widehat{{\rm d}^ne_\mu (0;z)}(\theta )=\left\langle P_n^\mu
(z),\theta ^{\otimes n}\right\rangle $ where $P_n^\mu (z)\in \left( {\cal N}%
^{\hat \otimes n}_{\Ckl }\right) ^{\prime }$. The kernel theorem really
gives a little more: $P_n^\mu (z)\in {\cal H}_{-p^{\prime }}^{\hat \otimes n}
$ for any $p^{\prime }(>p\geq p_0)$ such that the embedding operator $%
i_{p^{\prime },p}:{\cal H}_{p^{\prime }}\hookrightarrow {\cal H}_{p\text{ }}$%
is Hilbert-Schmidt. Thus we have 
\begin{equation}
\label{Pgenerator}e_\mu (\theta ;z)=\sum_{n=0}^\infty \frac
1{n!}\left\langle P_n^\mu (z),\theta ^{\otimes n}\right\rangle \quad \text{%
for }\theta \in {\cal U}_0,\ z\in {\cal N}_{\Ckl }^{\prime }\text{ .} 
\end{equation}
We will also use the notation%
$$
P_n^\mu (\varphi ^{(n)})(z):=\left\langle P_n^\mu (z),\varphi
^{(n)}\right\rangle ,\qquad \varphi ^{(n)}\in {\cal N}_{\Ckl }^{\hat \otimes
n},\quad n\in {\N}. 
$$
Thus for any measure satisfying Assumption 1 we have defined the ${\p}^\mu $%
-system

$$
{\p}^\mu =\left\{ \left\langle P_n^\mu (\cdot ),\varphi ^{(n)}\right\rangle
\ \bigg|\ \varphi ^{(n)}\in {\cal N}_{\Ckl }^{\hat \otimes n},\ n\in {\N}%
\right\} . 
$$
\bigskip\ 

Let us collect some properties of the polynomials $P_n^\mu (z).$

\begin{proposition}
For $x,y\in {\cal N}^{\prime }\ ,\ n\in \N $ the following holds\\[6mm]
(P1) \nopagebreak \hfill \\[-11mm] 
\begin{equation}
\label{(P1)}P_n^\mu (x)=\sum_{k=0}^n\binom nkx^{\otimes k}\hat \otimes
P_{n-k}^\mu (0),
\end{equation}
\\[5mm]
(P2)\nopagebreak \hfill \\[-11mm] 
\begin{equation}
\label{(P2)}x^{\otimes n}=\sum_{k=0}^n\binom nkP_k^\mu (x)\hat \otimes {\rm M%
}_{n-k}^\mu 
\end{equation}
\\[5mm]
(P3)\nopagebreak \hfill \\[-11mm] 
$$
P_n^\mu (x+y)=\sum_{k+l+m=n}\frac{n!}{k!\,l!\,m!}P_k^\mu (x)\hat \otimes
P_l^\mu (y)\hat \otimes {\rm M}_m^\mu  
$$
\begin{equation}
\label{(P3)}=\sum_{k=0}^n\binom nkP_k^\mu (x)\hat \otimes y^{\otimes (n-k)}
\end{equation}
(P4)\ Further we observe 
\begin{equation}
\label{(P4)}\E _\mu (\langle P_m^\mu (\cdot ),\varphi ^{(m)}\rangle
)=0\qquad \text{for }m\neq 0\ ,\varphi ^{(m)}\in {\cal N}_{\Ckl }^{\hat
\otimes m}\ .
\end{equation}
(P5) For all $p>p_0$ such that the embedding ${\cal H}_p\hookrightarrow 
{\cal H}_{p_0}$ is Hilbert--Schmidt and for all $\varepsilon >0$ small
enough $\left( \varepsilon \leq \frac{2^{-q_0}}{e\left\| i_{p,p_0}\right\|
_{HS}}\right) $ there exists a constant $C_{p,\varepsilon }>0$ with 
\begin{equation}
\label{(P5)}\left| P_n^\mu (z)\right| _{-p}\leq C_{p,\varepsilon
}\,n!\,\varepsilon ^{-n}\,e^{\varepsilon |z|_{-p}},\quad z\in {\cal H}_{-p,%
\Ckl }
\end{equation}
\end{proposition}

\TeXButton{Proof}{\proof}We restrict ourselves to a sketch of proof, details
can be found in \cite{ADKS94}.

\noindent (P1) This formula can be obtained simply by substituting 
\begin{equation}
\label{1/L}\frac 1{l_\mu (\theta )}=\sum\limits_{n=0}^\infty \frac
1{n!}\left\langle P_n^\mu (0),\theta ^{\otimes n}\right\rangle ,\quad \theta
\in {\cal N}_{\Ckl },\left| \theta \right| _q<\delta 
\end{equation}
and 
$$
e^{\left\langle x,\theta \right\rangle }=\sum\limits_{n=0}^\infty \frac
1{n!}\left\langle x^{\otimes n},\theta ^{\otimes n}\right\rangle ,\quad
\theta \in {\cal N}_{\Ckl },x\in {\cal N}^{\prime } 
$$
in the equality $e_\mu (\theta ;x)=e^{\left\langle x,\theta \right\rangle
}l_\mu ^{-1}(\theta )$. A comparison with (\ref{Pgenerator}) proves (P1).
The proof of (P2) is completely analogous to the proof of (P1).

\noindent (P3) We start from the following obvious equation of the
generating functions%
$$
e_\mu (\theta ;x+y)=e_\mu (\theta ;x)\,e_\mu (\theta ;y)\,l_\mu (\theta ) 
$$
This implies%
$$
\sum_{n=0}^\infty \frac 1{n!}\langle P_n^\mu (x+y),\theta ^{\otimes
n}\rangle =\sum_{k,l,m=0}^\infty \frac 1{k!\,l!\,m!}\,\langle P_k(x)\hat
\otimes P_l(y)\hat \otimes {\rm M}_m,\;\theta ^{\otimes (k+l+m)}\rangle 
$$
from this (P3) follows immediately.

\noindent (P4) To see this we use, $\theta \in {\cal N}_{\Ckl }$,%
$$
\sum_{n=0}^\infty \frac 1{n!}\E _\mu (\langle P_m^\mu (\cdot ),\theta
^{\otimes n}\rangle )=\E _\mu (e_\mu (\theta ;\cdot ))=\frac{\E _\mu
(e^{\langle \cdot ,\theta \rangle })}{l_\mu (\theta )}=1\ . 
$$
Then a comparison of coefficients and the polarization identity gives the
above result.

\noindent (P5) We can use 
\begin{equation}
\label{Pnxnorm}|P_n^\mu (z)|_{-p^{\prime }}\leq n!\left( \sup _{|\theta
|_p=\rho }\frac 1{l_\mu (\theta )}\right) e^{\rho |z|_{-p}}\left( \frac
e\rho \,\left\| i_{p^{\prime },p}\right\| _{HS}\right) ^n\;,\quad z\in {\cal %
H}_{-p,\Ckl } 
\end{equation}
$p>p_0,p^{\prime },\rho $ defined above. (\ref{Pnxnorm}) is a simple
consequence of the kernel theorem by (\ref{Cauchyemy}). In particular we
have 
$$
\left| P_n^\mu (0)\right| _{-p}\leq n!\left( \sup _{|\theta |_{p_0}=\rho
}\frac 1{l_\mu (\theta )}\right) \left( \frac e\rho \left\|
i_{p,p_0}\right\| _{HS}\right) ^n 
$$
If $p>p_0$ such that $\left\| i_{p,p_0}\right\| _{HS}$ is finite. For $%
0<\varepsilon \leq 2^{-q_0}/e\left\| i_{p,p_0}\right\| _{HS}$ we can fix $%
\rho =\varepsilon \,e\,\left\| i_{p,p_0}\right\| _{HS}\leq 2^{-q_0}$. With 
$$
C_{p,\varepsilon }:=\sup _{|\theta |_{p_0}=\rho }\frac 1{l_\mu (\theta )} 
$$
we have%
$$
\left| P_n^\mu (0)\right| _{-p}\leq C_{p,\varepsilon }\,n!\,\varepsilon
^{-n}. 
$$
Using (\ref{(P1)}) the following estimates hold 
\begin{eqnarray*}
\left| P_n^\mu (z)\right| _{-p} & \leq & \sum_{k=0}^n\binom nk\left| P_k^\mu %
(0)\right| _{-p}\left| z\right| _{-p}^{n-k}\ ,\qquad z\in {\cal H}_{-p,\Ckl } %
\\& \leq & C_{p,\varepsilon }\sum_{k=0}^n\tbinom nkk!\,\varepsilon ^{-k}\left| %
z\right| _{-p}^{n-k} %
\\&=& C_{p,\varepsilon }\,n!\,\varepsilon ^{-n}\sum_{k=0}^n\tfrac 1{(n-k)!}%
(\varepsilon \,|z|_{-p})^{n-k} %
\\& \leq & C_{p,\varepsilon }\,n!\,%
\varepsilon ^{-n}\,e^{\varepsilon |z|_{-p}}\ . 
\end{eqnarray*}
This completes the proof.\TeXButton{End Proof}{\endproof} \bigskip

\TeXButton{Note }{\note } The formulae (\ref{(P1)}) and (\ref{1/L}) can also be used as an
alternative definition of the polynomials $P_n^\mu (x)$\bigskip\ .

\example 
\newcounter{GaussA} \setcounter{GaussA}{\value{example}} \label{GaussAP} Let
us compare to the case of Gaussian Analysis. Here one has 
$$
l_{\gamma _{{\cal H}}}(\theta )=e^{\frac 12\langle \theta ,\theta \rangle }\
,\qquad \theta \in {\cal N}_{\Ckl } 
$$
Then it follows%
$$
{\rm M}_{2n}^\mu =(-1)^nP_{2n}^\mu (0)=\frac{(2n)!}{n!\,2^n}{\rm Tr}^{\hat
\otimes n}\ ,\qquad n\in \N  
$$
and ${\rm M}_n^\mu =P_n^\mu (0)=0$ if $n$ is odd. Here ${\rm Tr}\in {\cal N}%
^{\prime \otimes 2}$ denotes the trace kernel defined by 
\begin{equation}
\label{Trace}\langle {\rm Tr},\eta \otimes \xi \rangle =(\eta ,\xi )\
,\qquad \eta ,\xi \in {\cal N} 
\end{equation}
A simple comparison shows that 
$$
P_n^\mu (x)=:x^{\otimes n}:\text{ } 
$$
and%
$$
e_\mu (\theta ;x)=:e^{\langle x,\theta \rangle }: 
$$
where the r.h.s. denotes usual Wick ordering see e.g., \cite{BeKo88,HKPS93}.
This procedure is uniquely defined by%
$$
\langle :x^{\otimes n}:,\xi ^{\otimes n}\rangle =2^{-\frac n2}|\xi
|^n\,H_n\left( \tfrac 1{\sqrt{2}|\xi |}\langle x,\xi \rangle \right) \
,\qquad \xi \in {\cal N} 
$$
where $H_n$ denotes the Hermite polynomial of order $n$ (see e.g., \cite
{HKPS93} for the normalization we use).\bigskip\ 

Now we are ready to give the announced description of ${\cal P}({\cal N}%
^{\prime })$.

\begin{lemma}
\label{PrepLemma}For any $\varphi \in {\cal P}({\cal N}^{\prime })$ there
exists a unique representation 
\begin{equation}
\label{Prep}\varphi (x)=\sum\limits_{n=0}^N\left\langle P_n^\mu (x),\varphi
^{(n)}\right\rangle \ ,\quad \text{ }\varphi ^{(n)}\in {\cal N}_{\Ckl %
}^{\hat \otimes n}
\end{equation}
and vice versa, any functional of the form (\ref{Prep}) is a smooth
polynomial.
\end{lemma}

\TeXButton{Proof}{\proof}The representations from Definition \ref{P(N)} and
equation (\ref{Prep}) can be transformed into one another using (\ref{(P1)})
and (\ref{(P2)}). \TeXButton{End Proof}{\endproof}

\subsection{The dual Appell system and the representation theorem for 
\texorpdfstring{${\cal P}_\mu ^{\prime }({\cal N}^{\prime })$} {P\textacute (N\textacute )  }
}

To give an internal description of the type (\ref{Prep}) for ${\cal P}_\mu
^{\prime }({\cal N}^{\prime })$ we have to construct an appropriate system
of generalized functions, the ${\Q}^\mu $-system. The construction we
propose here is different from that of \cite{ADKS94} where smoothness of the
logarithmic derivative of $\mu $ was demanded and used for the construction
of the ${\Q}^\mu $-system. To avoid this additional assumption (which
excludes e.g., Poisson measures) we propose to construct the ${\Q}^\mu $%
-system using differential operators.

Define a differential operator of order $n$ with constant coefficient $\Phi
^{(n)}\in \left( {\cal N}_{\Ckl }^{\hat \otimes n}\right) ^{\prime }$ 
$$
D(\Phi ^{(n)})\langle x^{\otimes m},\varphi ^{(m)}\rangle =\QDATOPD\{ . {%
\dfrac{m!}{(m-n)!}\langle x^{\otimes (m-n)}\hat \otimes \Phi ^{(n)},\varphi
^{(m)}\rangle \text{ \quad for }m\geq n}{\text{\hspace{1.5cm}}0\text{%
\hspace{4cm}for }m<n} 
$$
($\varphi ^{(m)}\in {\cal N}_{\Ckl }^{\hat \otimes m},m\in {\N}$) and extend
by linearity from the monomials to ${\cal P}({\cal N}^{\prime }).$

\begin{lemma}
\label{Dcontinuity}$D(\Phi ^{(n)})$is a continuous linear operator from $%
{\cal P}({\cal N}^{\prime })$ to ${\cal P}({\cal N}^{\prime })$ .
\end{lemma}

\TeXButton{Remark }{\remark } For $\Phi ^{(1)}\in {\cal N}^{\prime }$ we
have the usual G\^ateaux derivative as e.g., in white noise analysis \cite
{HKPS93}%
$$
D(\Phi ^{(1)})\varphi =D_{\Phi ^{(1)}}\varphi :=\frac{{\rm d}}{{\rm d}t}%
\varphi (\cdot +t\Phi ^{(1)})|_{t=0} 
$$
for $\varphi \in {\cal P}({\cal N})$ and we have $D(\left( \Phi
^{(1)}\right) ^{\otimes n})=(D_{\Phi ^{(1)}})^n$ thus $D(\left( \Phi
^{(1)}\right) ^{\otimes n})$ is in fact a differential operator of order $n$.

\TeXButton{Proof}{\proof}By definition ${\cal P}({\cal N}^{\prime })$ is
isomorphic to the topological direct sum of tensor powers ${\cal N}_{\Ckl %
}^{\hat \otimes n}$%
$$
{\cal P}({\cal N}^{\prime })\simeq \bigoplus_{n=0}^\infty {\cal N}_{\Ckl %
}^{\hat \otimes n}\text{ .} 
$$
Via this isomorphism $D(\Phi ^{(n)})$ transforms each component ${\cal N}_{%
\Ckl }^{\hat \otimes m}$, $m\geq n$ by 
$$
\varphi ^{(m)}\mapsto \frac{n!}{(m-n)!}(\Phi ^{(n)},\;\varphi ^{(m)})_{{\cal %
H}^{\hat \otimes n}} 
$$
where the contraction $(\Phi ^{(n)},\;\varphi ^{(m)})_{{\cal H}^{\hat
\otimes n}}$ $\in {\cal N}_{\Ckl }^{\otimes (m-n)}$ is defined by 
\begin{equation}
\label{contraction}\langle x^{\otimes (m-n)},\;(\Phi ^{(n)},\;\varphi
^{(m)})_{{\cal H}^{\hat \otimes n}}\rangle :=\langle x^{\otimes (m-n)}\hat
\otimes \Phi ^{(n)},\varphi ^{(m)}\rangle 
\end{equation}
for all $x\in {\cal N}^{\prime }$. It is easy to verify that%
$$
|(\Phi ^{(n)},\;\varphi ^{(m)})_{{\cal H}^{\hat \otimes n}}|_q\leq |\Phi
^{(n)}|_{-q}|\varphi ^{(m)}|_q\text{ ,\qquad }q\in \N  
$$
which guarantees that $(\Phi ^{(n)},\;\varphi ^{(m)})_{{\cal H}^{\hat
\otimes n}}\in {\cal N}_{\Ckl }^{\otimes (m-n)}$ and shows at the same time
that $D(\Phi ^{(n)})$ is continuous on each component. This is sufficient to
ensure the stated continuity of $D(\Phi ^{(n)})$ on ${\cal P}({\cal N}%
^{\prime }).$\TeXButton{End Proof}{\endproof}\bigskip\ 

\begin{lemma}
For $\Phi ^{(n)}\in {\cal N}_{\Ckl }^{\prime \hat \otimes n}$ , $\varphi
^{(m)}\in {\cal N}_{\Ckl }^{\hat \otimes m}$ we have\\[9mm](P6)\\[-14mm] 
\begin{equation}
\label{(P6)}D(\Phi ^{(n)})\langle P_m^\mu (x),\varphi ^{(m)}\rangle
=\QATOPD\{ . {\dfrac{m!}{(m-n)!}\left\langle P_{m-n}^\mu (x)\hat \otimes
\Phi ^{(n)},\;\varphi ^{(m)}\right\rangle \text{ for }m\geq n}{\text{%
\hspace{1.5cm}}0\text{\hspace{4.1cm} for }m<n}
\end{equation}
\end{lemma}

\TeXButton{Proof}{\proof}This follows from the general property of Appell
polynomials which behave like ordinary powers under differentiation. More
precisely, by using 
$$
\langle P_m^\mu ,\theta ^{\otimes m}\rangle =\left. \left( \frac{{\rm d}}{%
{\rm d}t}\right) ^me_\mu (t\theta ;\cdot )\right| _{t=0}\ ,\qquad \theta \in 
{\cal N}_{\Ckl } 
$$
we have%
\begin{eqnarray*}
D(\Phi ^{(1)})\langle P_m^\mu (x),\theta ^{\otimes m}\rangle %
&=& \left. \frac{{\rm d}}{{\rm d}\lambda }\langle P_m^\mu (x+\lambda %
\Phi ^{(1)}),\theta ^{\otimes m}\rangle \right| _{\lambda =0} %
\\&=& \left. \left( \frac \partial {\partial t}\right) ^m\frac \partial {\partial \lambda }e_\mu (t\theta ;x+\lambda \Phi ^{(1)})%
\right| _{\T \QATOP{t=0}{\lambda =0}} %
\\&=& \langle \Phi ^{(1)},\theta \rangle \left. \left( \tfrac \partial {\partial t}\right) ^mt\;e_\mu (t\theta ;x)\right| _{t=0} %
\\&=& \left. \langle \Phi ^{(1)},\theta \rangle \sum_{k=0}^m\tbinom mk\left( %
\left( \tfrac{{\rm d}}{{\rm d}t}\right) ^kt\right) \left( \tfrac{{\rm d}}{{\rm d}t}\right) ^{m-k}e_\mu (t\theta ;x)\right| _{t=0} %
\\&=& \left. m\,\langle \Phi ^{(1)},\theta \rangle \left( \tfrac{{\rm d}}{{\rm d}t}\right) ^{m-1}e_\mu (t\theta ;x)\right| _{t=0} %
\\&=& m\,\langle \Phi ^{(1)},\theta \rangle \left\langle P_{m-1}^\mu (x),\theta %
^{\otimes (m-1)}\right\rangle \text{ .} 
\end{eqnarray*}
This proves 
$$
D(\Phi ^{(1)})\langle P_m^\mu ,\varphi ^{(m)}\rangle =m\left\langle
P_{m-1}^\mu \hat \otimes \Phi ^{(1)},\;\varphi ^{(m)}\right\rangle . 
$$
The property (\ref{(P6)}), then follows by induction.\TeXButton{End Proof}
{\endproof}\bigskip\ 

In view of Lemma \pageref{Dcontinuity} it is possible to define the adjoint
operator $D(\Phi ^{(n)})^{*}:{\cal P}_\mu ^{\prime }({\cal N}^{\prime
})\rightarrow {\cal P}_\mu ^{\prime }({\cal N}^{\prime })$ for $\Phi
^{(n)}\in {\cal N}_{\Ckl }^{\prime \hat \otimes n}$ . Further we can
introduce the constant function $\1 \in {\cal P}_\mu ^{\prime }({\cal N}%
^{\prime })$ such that $\1 (x)\equiv 1$ for all $x\in {\cal N}^{\prime }$ ,
so 
$$
\langle \!\langle \1 ,\;\varphi \rangle \!\rangle _\mu =\int_{{\cal N}%
^{\prime }}\varphi (x)\,{\rm d\mu }(x)=\E _\mu (\varphi ). 
$$
Now we are ready to define our $\Q$-system.

\begin{definition}
For any $\Phi ^{(n)}\in \left( {\cal N}_{\Ckl }^{\hat \otimes n}\right)
^{\prime }$ we define $Q_n^\mu (\Phi ^{(n)})\in {\cal P}_\mu ^{\prime }(%
{\cal N}^{\prime })$ by%
$$
Q_n^\mu (\Phi ^{(n)})=D(\Phi ^{(n)})^{*}{\1}\ . 
$$
\end{definition}

We want to introduce an additional formal notation $Q_n^\mu (x)$ which
stresses the linearity of $\Phi ^{(n)}\mapsto Q_n^\mu (\Phi ^{(n)})\in P_\mu
^{\prime }({\cal N}^{\prime }):$%
$$
\langle Q_n^\mu ,\Phi ^{(n)}\rangle :=Q_n^\mu (\Phi ^{(n)})\ . 
$$

\example It is possible to put further assumptions on the measure $\mu $ to
ensure that the expression is more than formal. Let us assume a smooth
measure (i.e., the logarithmic derivative of $\mu $ is infinitely
differentiable, see \cite{ADKS94} for details) with the property%
$$
\exists q\in \N \ ,\ \exists \{C_n\geq 0,\;n\in \N \}:\forall \xi \in {\cal N%
} 
$$
$$
\left| \int D_\xi ^n\varphi \;{\rm d}\mu (x)\right| \leq C_n\left\| \varphi
\right\| _{L^2(\mu )}|\xi |_q^n 
$$
where $\varphi $ is any finitely based bounded ${\cal C}^\infty $-function
on ${\cal N}^{\prime }$. This obviously establishes a bound of the type 
$$
\left\| Q_n^\mu (\xi _1\otimes \cdots \otimes \xi _n)\right\| _{L^2(\mu
)}\leq C_n^{\prime }\prod_{j=1}^n|\xi _j|_q\ ,\qquad \xi _1,\ldots ,\xi
_n\in {\cal N\ },\ n\in \N  
$$
which is sufficient to show (by means of kernel theorem) that there exists $%
Q_n^\mu (x)\in \left( {\cal N}_{\Ckl }^{\hat \otimes n}\right) ^{\prime }$
for almost all $x\in {\cal N}^{\prime }$ such that we have the representation%
$$
Q_n^\mu (\varphi ^{(n)})(x)=\langle Q_n^\mu (x),\varphi ^{(n)}\rangle \
,\qquad \varphi ^{(n)}\in {\cal N}_{\Ckl }^{\hat \otimes n} 
$$
for almost all $x\in {\cal N}^{\prime }$. For any smooth kernel $\varphi
^{(n)}\in {\cal N}_{\Ckl }^{\hat \otimes n}$ we have then that the function%
$$
x\mapsto \langle Q_n^\mu (x),\varphi ^{(n)}\rangle \ =Q_n^\mu \left( \varphi
^{(n)}\right) (x)\ 
$$
belongs to $L^2(\mu ).$\bigskip\ 

\example The simplest non trivial case can be studied using finite
dimensional real analysis. We consider $\R $ as our basic Hilbert space and
as our nuclear space ${\cal N}$. Thus the nuclear ``triple'' is simply%
$$
\R \subseteq \R \subseteq \R  
$$
and the dual pairing between a ``test function'' and a ``distribution''
degenerates to multiplication. On $\R $ we consider a measure ${\rm d}\mu
(x)=\rho (x)\,{\rm d}x$ where $\rho $ is a positive ${\cal C}^\infty $%
--function on $\R $ such that Assumptions 1 and 2 are fulfilled. In this
setting the adjoint of the differentiation operator is given by%
$$
\left( \frac{{\rm d}}{{\rm d}x}\right) ^{*}f(x)=-\left( \tfrac{{\rm d}}{{\rm %
d}x}+\beta (x)\right) f(x)\ ,\qquad f\in {\cal C}^1(\R ) 
$$
where the logarithmic derivative $\beta $ of the measure $\mu $ is given by 
$$
\beta =\frac{\rho ^{\prime }}\rho 
$$
This enables us to calculate the $\Q ^\mu $-system. One has 
$$
Q_n^\mu (x)=\left( \left( \tfrac{{\rm d}}{{\rm d}x}\right) ^{*}\right) ^n\1 %
=(-1)^n\left( \tfrac{{\rm d}}{{\rm d}x}+\beta (x)\right) ^n\1  
$$
$$
=(-1)^n\frac{\rho ^{(n)}(x)}{\rho (x)}\ . 
$$
The last equality can be seen by simple induction. \\ If $\rho =\frac 1{%
\sqrt{2\pi }}e^{-\frac 12x^2}$ is the Gaussian density $Q_n^\mu $ is related
to the n$^{th}$ Hermite polynomial:%
$$
Q_n^\mu (x)=2^{-n/2}H_n\left( \tfrac x{\sqrt{2}}\right) \ . 
$$
\bigskip\ 

\begin{definition}
We define the $\Q^\mu $-system in ${\cal P}_\mu ^{\prime }({\cal N}^{\prime
})$ by 
$$
\Q ^\mu =\left\{ Q_n^\mu (\Phi ^{(n)})\ \Big|\ \qquad \Phi ^{(n)}\in \left( 
{\cal N}_{\Ckl }^{\hat \otimes n}\right) ^{\prime },\ n\in \N _0\ \right\} \
, 
$$
and the pair $(\p ^\mu ,\Q^\mu )$ will be called the Appell system $\A ^\mu $
generated by the measure $\mu $.
\end{definition}

Now we are going to discuss the central property of the Appell system $\A %
^\mu $.

\begin{theorem}
\label{BiorTh}{\rm (Biorthogonality w.r.t. }$\mu ${\rm )} 
\begin{equation}
\label{QnPnPair}\left\langle \!\!\left\langle \langle Q_n^\mu (\Phi
^{(n)}),\ \langle P_m^\mu ,\varphi ^{(m)}\rangle \right\rangle
\!\!\right\rangle _\mu =\delta _{m,n}\;n!\;\langle \Phi ^{(n)},\varphi
^{(n)}\rangle 
\end{equation}
for $\Phi ^{(n)}\in \left( {\cal N}_{\Ckl }^{\hat \otimes n}\right) ^{\prime
}$ and $\varphi ^{(m)}\in {\cal N}_{\Ckl }^{\hat \otimes m}$ .
\end{theorem}

\TeXButton{Proof}{\proof}It follows from (\ref{(P4)}) and (\ref{(P6)}) that%
\begin{eqnarray*}
\left\langle \!\!\left\langle Q_n^\mu (\Phi ^{(n)}),\ \langle P_m^\mu %
,\varphi ^{(m)}\rangle \right\rangle \!\!\right\rangle _\mu %
&=& \left\langle \!\!\left\langle \1 ,D(\Phi ^{(n)})\langle P_m^\mu ,\varphi %
^{(m)}\rangle \right\rangle \!\!\right\rangle _\mu %
\\&=& \frac{m!}{(m-n)!}\E _\mu \left( \langle P_{(m-n)}^\mu \hat \otimes \Phi %
^{(n)},\;\varphi ^{(m)}\rangle \right) %
\\&=& m!\;\delta _{m,n}\;\langle \Phi ^{(m)},\varphi ^{(m)}\rangle \ . 
\end{eqnarray*}
\TeXButton{End Proof}{\endproof}\bigskip\ 

Now we are going to characterize the space ${\cal P}_\mu ^{\prime }({\cal N}%
^{\prime })$

\begin{theorem}
\label{PStrichRep}For all $\Phi \in {\cal P}_\mu ^{\prime }({\cal N}^{\prime
})$ there exists a unique sequence $\{\Phi ^{(n)}\big|\ n\in \N _0\},\ \Phi
^{(n)}\in \left( {\cal N}_{\Ckl }^{\hat \otimes n}\right) ^{\prime }$ such
that 
\begin{equation}
\label{Qexpansion}\Phi =\sum_{n=0}^\infty Q_n^\mu (\Phi ^{(n)})\equiv
\sum_{n=0}^\infty \langle Q_n^\mu ,\Phi ^{(n)}\rangle 
\end{equation}
and vice versa, every series of the form (\ref{Qexpansion}) generates a
generalized function in ${\cal P}_\mu ^{\prime }({\cal N}^{\prime }).$
\end{theorem}

\TeXButton{Proof}{\proof}For $\Phi \in {\cal P}_\mu ^{\prime }({\cal N}%
^{\prime })$ we can uniquely define $\Phi ^{(n)}\in \left( {\cal N}_{\Ckl %
}^{\hat \otimes n}\right) ^{\prime }$ by 
$$
\langle \Phi ^{(n)},\varphi ^{(n)}\rangle =\frac 1{n!}\ \!\langle \!\langle
\Phi ,\;\langle P_n^\mu ,\varphi ^{(n)}\rangle \rangle \!\rangle _\mu \
,\qquad \varphi ^{(n)}\in {\cal N}_{\Ckl }^{\hat \otimes n} 
$$
This definition is possible because $\langle P_n^\mu ,\varphi ^{(n)}\rangle
\in {\cal P}({\cal N}^{\prime })$. The continuity of $\varphi ^{(n)}\mapsto
\langle \Phi ^{(n)},\varphi ^{(n)}\rangle $ follows from the continuity of $%
\varphi \mapsto \langle \!\langle \Phi ,\varphi \rangle \!\rangle \ ,\
\varphi \in $ ${\cal P}({\cal N}^{\prime })$. This implies that $\varphi
\mapsto \sum_{n=0}^\infty n!\;\langle \Phi ^{(n)},\varphi ^{(n)}\rangle $ is
continuous on ${\cal P}({\cal N}^{\prime })$. This defines a generalized
function in ${\cal P}_\mu ^{\prime }({\cal N}^{\prime })$, which we denote
by $\sum_{n=0}^\infty Q_n^\mu (\Phi ^{(n)})$. In view of Theorem \ref{BiorTh}
it is obvious that 
$$
\Phi =\sum_{n=0}^\infty Q_n^\mu (\Phi ^{(n)})\ . 
$$

To see the converse consider a series of the form (\ref{Qexpansion}) and $%
\varphi \in {\cal P}({\cal N}^{\prime })$. Then there exist $\varphi
^{(n)}\in {\cal N}_{\Ckl }^{\hat \otimes n}\ ,\ n\in \N $ and $N\in \N $
such that we have the representation%
$$
\varphi =\sum_{n=0}^NP_n^\mu (\varphi ^{(n)})\ . 
$$
So we have%
$$
\!\left\langle \!\!\!\left\langle \sum_{n=0}^\infty Q_n^\mu (\Phi
^{(n)}),\varphi \right\rangle \!\!\!\right\rangle _{\!\!\mu
}:=\sum_{n=0}^Nn!\;\langle \Phi ^{(n)},\varphi ^{(n)}\rangle 
$$
because of Theorem \ref{BiorTh}. The continuity of $\varphi \mapsto \langle
\!\langle \sum_{n=0}^\infty Q_n^\mu (\Phi ^{(n)}),\varphi \rangle \!\rangle
_\mu $ follows because $\varphi ^{(n)}\mapsto \langle \Phi ^{(n)},\varphi
^{(n)}\rangle $ is continuous for all $n\in \N $ .\TeXButton{End Proof}
{\endproof}

%% file: alms.tex
\LaTeXparent{nga4.tex}

\section{Test functions on a linear space with measure \label{testfunctions}}

In this section we will construct the test function space $({\cal N})^1$ and
study its properties. On the space ${\cal P}({\cal N}^{\prime })$ we can
define a system of norms using the representation from Lemma \ref{PrepLemma}%
. Let 
$$
\varphi =\sum_{n=0}^N\langle P_n^\mu ,\;\varphi ^{(n)}\rangle \in {\cal P}(%
{\cal N}^{\prime }) 
$$
be given, then $\varphi ^{(n)}\in {\cal H}_{p,\Ckl }^{\hat \otimes n}$ for
each $p\geq 0\ \ (n\in \N )$. Thus we may define for any $p,q\in \N $ a
Hilbertian norm on ${\cal P}({\cal N}^{\prime })$ by 
$$
\left\| \varphi \right\| _{p,q,\mu }^2=\sum_{n=0}^N(n!)^2\;2^{nq}\;|\varphi
^{(n)}|_p^2 
$$
The completion of ${\cal P}({\cal N}^{\prime })$ w.r.t. $\left\| \cdot
\right\| _{p,q,\mu }$ is called $({\cal H}_p)_{q,\mu }^1$ .

\begin{definition}
We define 
$$
({\cal N})_\mu ^1:=\ \stackunder{p,q\in \N }{\rm pr\ lim}({\cal H}_p)_{q,\mu
}^1\ . 
$$
\end{definition}

This space has the following properties

\begin{theorem}
$({\cal N})_\mu ^1$ is a nuclear space. The topology $({\cal N})_\mu ^1$ is
uniquely defined by the topology on ${\cal N}$: It does not depend on the
choice of the family of norms $\{|\cdot |_p\}$.
\end{theorem}

\TeXButton{Proof}{\proof}Nuclearity of $({\cal N})_\mu ^1$ follows
essentially from that of ${\cal N}.$ For fixed $p,q$ consider the embedding 
$$
I_{p^{\prime },q^{\prime },p,q}:\left( {\cal H}_{p^{\prime }}\right)
_{q^{\prime },\mu }^1\rightarrow \left( {\cal H}_p\right) _{q,\mu }^1 
$$
where $p^{\prime }$ is chosen such that the embedding%
$$
i_{p^{\prime },p}:\,{\cal H}_{p^{\prime }}\rightarrow {\cal H}_p 
$$
is Hilbert--Schmidt. Then $I_{p^{\prime },q^{\prime },p,q}$ is induced by%
$$
I_{p^{\prime },q^{\prime },p,q}\varphi =\sum_{n=0}^\infty \langle P_n^\mu
,i_{p^{\prime },p}^{\otimes n}\varphi ^{(n)}\rangle \quad \text{ for \quad }%
\varphi =\sum_{n=0}^\infty \langle P_n^\mu ,\varphi ^{(n)}\rangle \in \left( 
{\cal H}_{p^{\prime }}\right) _{q^{\prime },\mu }^1. 
$$
Its Hilbert--Schmidt norm is easily estimated by using an orthonormal basis
of $\left( {\cal H}_{p^{\prime }}\right) _{q^{\prime },\mu }^1$. The result
is the bound 
$$
\left\| I_{p^{\prime },q^{\prime },p,q}\right\| _{HS}^2\le \sum_{n=0}^\infty
2^{n(q-q^{\prime })}\left\| i_{p^{\prime },p}\right\| _{HS}^{2n} 
$$
which is finite for suitably chosen $q^{\prime }$.

Let us assume that we are given two different systems of Hilbertian norms $%
\left| \,\cdot \,\right| _p$ and $\left| \,\cdot \,\right| _k^{\prime }$,
such that they induce the same topology on ${\cal N}$ . For fixed $k$ and $l$
we have to estimate $\left\| \,\cdot \,\right\| _{k,l,\mu }^{\prime }$ by $%
\left\| \,\cdot \,\right\| _{p,q,\mu }$ for some $p,q$ (and vice versa which
is completely analogous). Since $\left| \,\cdot \,\right| _k^{\prime }$ has
to be continuous with respect to the projective limit topology on ${\cal N}$%
, there exists $p$ and a constant $C$ such that $\left| f\right| _k^{\prime
}\leq C\left| f\right| _p$, for all $f\in {\cal N}$, i.e., the injection $i$
from ${\cal H}_p$ into the completion ${\cal K}_k$ of ${\cal N}$ with
respect to $|\,\cdot \,|_k^{\prime }$ is a mapping bounded by $C$. We denote
by $i$ also its linear extension from ${\cal H}_{p,{\,}\Ckl }$ into ${\cal K}%
_{{\,}k,\Ckl }$. It follows that $i^{\otimes n}$ is bounded by $C^n$ from $%
{\cal H}_{{\,}p,\Ckl }^{\otimes n}$ into ${\cal K}_{{\,}k,\Ckl   }^{\otimes
n}$. Now we choose $q$ such that $2^{{\frac{q-l}2}}\geq C$. Then

\begin{eqnarray*}
\left\| \,\cdot \,\right\| _{k,l,\mu }^{\prime 2}&=&\sum_{n=0}^\infty %
(n!)^2\,2^{nl}\left| \,\cdot \,\right| _k^{\prime 2} %
\\&\leq & \sum_{n=0}^\infty (n!)^2\,2^{nl}C^{2n}\left| \,\cdot \,%
\right| _p^2 %
\\&\leq &\left\| \,\cdot \,\right\| _{p,q,\mu }^2\ , 
\end{eqnarray*}
which had to be proved.\TeXButton{End Proof}{\endproof}\bigskip\ 

\begin{lemma}
\label{L2NormPn}There exist $p,C,K>0$ such that for all $n$ 
\begin{equation}
\label{Pnx2norm}\int |P_n^\mu (x)|_{-p}^2\;{\rm d}\mu (x)\leq (n!)^2\,C^n\,K
\end{equation}
\end{lemma}

\TeXButton{Proof}{\proof} The estimate (\ref{Pnxnorm}) may be used for $\rho
\leq 2^{-q_0}$ and $\rho \leq 2\varepsilon _\mu $ ($\varepsilon _\mu $ from
Lemma \ref{equiLemma}).\\This gives 
$$
\int |P_n^\mu (x)|_{-p}^2\;{\rm d}\mu (x)\leq (n!)^2\left( \frac e\rho
\left\| i_{p,p_0}\right\| _{HS}\right) ^{2n}\int e^{2\rho |x|_{-p_0}}{\rm d}%
\mu (x) 
$$
which is finite because of Lemma \ref{equiLemma}.\TeXButton{End Proof}
{\endproof}

\begin{theorem}
There exist $p^{\prime },q^{\prime }>0$ such that for all $p\geq p^{\prime
},\ q\geq q^{\prime }$ the topological embedding $({\cal H}_p)_{q,\mu
}^1\subset L^2(\mu )$ holds.
\end{theorem}

\TeXButton{Proof}{\proof}Elements of the space $({\cal N})_\mu ^1$ are
defined as series convergent in the given topology. Now we need to study the
convergence of these series in $L^2(\mu )$. Choose $q^{\prime }$ such that $%
C>2^{q^{\prime }\text{ }}$ ($C$ from estimate (\ref{Pnx2norm})). Let us take
an arbitrary 
$$
\varphi =\sum_{n=0}^\infty \langle P_n^\mu ,\varphi ^{(n)}\rangle \in {\cal P%
}({\cal N}^{\prime }) 
$$
For $p>p^{\prime }$ ($p^{\prime }$ as in Lemma \ref{L2NormPn} ) and $%
q>q^{\prime }$ the following estimates hold%
\begin{eqnarray*}
\left\| \varphi \right\| _{L^2(\mu )}%
& \leq & \sum_{n=0}^\infty \left\| \langle%
P_n^\mu ,\varphi ^{(n)}\rangle \right\| _{L^2(\mu )}%
\\&\leq & \sum_{n=0}^\infty |\varphi ^{(n)}|_{-p}\left\| \,|P_n^\mu%
|_{-p}\right\| _{L^2(\mu )} %
\\&\leq & K\sum_{n=0}^\infty n!\,2^{nq/2}\left| \varphi ^{(n)}\right|%
_{-p}(C2^{-q})^{n/2} %
\\&\leq & K\left( \sum_{n=0}^\infty (C\,2^{-q})^n\right) ^{\frac 12}\left(%
\sum_{n=0}^\infty (n!)^2\,2^{qn}\left| \varphi ^{(n)}\right| _{-p}^2\right)%
^{\frac 12} %
\\&=& K\left( 1-C\,2^{-q}\right) ^{-1/2}\left\| \varphi \right\| _{p,q,\mu }%
\text{.}%
\end{eqnarray*}
Taking the closure the inequality extends to the whole space $({\cal H}%
_p)_q^1$.\TeXButton{End Proof}{\endproof}

\begin{corollary}
\label{N1inL2}$({\cal N})_\mu ^1$ is continuously and densely embedded in $%
L^2(\mu )$.
\end{corollary}

\example \newcounter{myexponent} \setcounter{myexponent}{\value{example}} 
{\it ($\mu $-exponentials as test functions)} \smallskip 
\\The $\mu $-exponential given in (\ref{Pgenerator}) has the following norm%
$$
||e_\mu (\theta ;\cdot )||_{p,q,\mu }^2=\sum_{n=0}^\infty 2^{nq}\,|\theta
|_p^{2n}\ ,\qquad \theta \in {\cal N}_{\Ckl } 
$$
This expression is finite if and only if $2^q|\theta |_p^2<1$. Thus we have $%
e_\mu (\theta ;\cdot )\notin ({\cal N})_\mu ^1$ if $\theta \neq 0$. But we
have that $e_\mu (\theta ;\cdot )$ is a test function of finite order i.e., $%
e_\mu (\theta ;\cdot )\in ({\cal H}_p)_q^1$ if $2^q|\theta |_p^2<1$. This is
in contrast to some useful spaces of test functions in Gaussian Analysis,
see e.g., \cite{BeKo88,HKPS93}.

The set of all $\mu $--exponentials $\{e_\mu (\theta ;\cdot )\;|\;2^q|\theta
|_p^2<1,\ \theta \in {\cal N}_{\Ckl }\}$ is a total set in $({\cal H}_p)_q^1$%
. This can been shown using the relation ${\rm d}^ne_\mu (0;\cdot )(\theta
_1,...,\theta _n)=\langle P_n^\mu ,\theta _1\hat \otimes \cdots \hat \otimes
\theta _n\rangle .$

\begin{proposition}
\label{N1inEmin}Any test function $\varphi $ in $({\cal N})_\mu ^1$ has a
uniquely defined extension to ${\cal N}_{\Ckl }^{\prime }$ as an element of $%
{\cal E}_{\min }^1\left( {\cal N}_{\Ckl }^{\prime }\right) $
\end{proposition}

\TeXButton{Proof}{\proof}Any element $\varphi $ in $({\cal N})_\mu ^1$ is
defined as a series of the following type 
$$
\varphi =\sum_{n=0}^\infty \langle P_n^\mu ,\varphi ^{(n)}\rangle \ ,\qquad
\varphi ^{(n)}\in {\cal N}_{\Ckl }^{\hat \otimes n} 
$$
such that%
$$
\left\| \varphi \right\| _{p,q,\mu }^2=\sum_{n=0}^\infty
(n!)^2\,2^{nq}\,|\varphi ^{(n)}|_p^2 
$$
is finite for each $p,q\in \N $ . In this proof we will show the convergence
of the series%
$$
\sum_{n=0}^\infty \langle P_n^\mu (z),\varphi ^{(n)}\rangle ,\quad z\in 
{\cal H}_{-p,\Ckl } 
$$
to an entire function in $z$.

Let $p>p_0$ such that the embedding $i_{p,p_0}:{\cal H}_p\hookrightarrow 
{\cal H}_{p_0}$ is Hilbert-Schmidt. Then for all $0<\varepsilon \leq
2^{-q_0}/e\left\| i_{p,p_0}\right\| _{HS}$ we can use (\ref{(P5)}) and
estimate as follows%
\begin{eqnarray*}
\sum_{n=0}^\infty |\langle P_n^\mu (z),\varphi ^{(n)}\rangle |%
& \leq & \sum_{n=0}^\infty |P_n^\mu (z)|_{-p}|\varphi ^{(n)}|_p %
\\& \leq & C_{p,\varepsilon }\,e^{\varepsilon |z|_{-p}}\sum_{n=0}^\infty%
n!\,|\varphi ^{(n)}|_p\,\varepsilon ^{-n} %
\\& \leq & C_{p,\varepsilon }\,\,e^{\varepsilon |z|_{-p}}\,\left(%
\sum_{n=0}^\infty (n!)^22^{nq}|\varphi ^{(n)}|_p^2\right) ^{1/2}\left(%
\sum_{n=0}^\infty 2^{-nq}\varepsilon ^{-2n}\right) ^{1/2} %
\\&=& C_{p,\varepsilon }\,\,\left( 1-2^{-q}\varepsilon ^{-2}\right)%
^{-1/2}\,\left\| \varphi \right\| _{p,q,\mu }\,\;e^{\varepsilon |z|_{-p}} %
\end{eqnarray*}
if $2^q>\varepsilon ^{-2}$. That means the series $\sum_{n=0}^\infty \langle
P_n^\mu (z),\varphi ^{(n)}\rangle $ converges uniformly and absolutely in
any neighborhood of zero of any space ${\cal H}_{-p,\Ckl }$ . Since each
term $\langle P_n^\mu (z),\varphi ^{(n)}\rangle $ is entire in $z$ the
uniform convergence implies that $z\mapsto \sum_{n=0}^\infty \langle P_n^\mu
(z),\varphi ^{(n)}\rangle $ is entire on each ${\cal H}_{-p,\Ckl }$ and
hence on ${\cal N}_{\Ckl }^{\prime }$. This completes the proof.%
\TeXButton{End Proof}{\endproof}\bigskip\ 

The following corollary is an immediate consequence of the above proof and
gives an explicit estimate on the growth of the test functions.

\begin{corollary}
\label{phi(z)Betrag}For all $p>p_0$ such that the norm $\left\|
i_{p,p_0}\right\| _{HS}$ of the embedding is finite and for all $%
0<\varepsilon \leq 2^{-q_0}/e\left\| i_{p,p_0}\right\| _{HS}$ we can choose $%
q\in \N $ such that $2^q>\varepsilon ^{-2}$ to obtain the following bound.%
$$
\left| \varphi (z)\right| \leq C\,\left\| \varphi \right\| _{p,q,\mu
}\,e^{\varepsilon |z|_{-p}}\ ,\qquad \varphi \in ({\cal N})_\mu ^1,\ z\in 
{\cal H}_{-p,\Ckl }\text{ ,} 
$$
where%
$$
C=C_{p,\varepsilon }\,\left( 1-2^{-q}\varepsilon ^{-2}\right) ^{-1/2}. 
$$
\end{corollary}

$\ $Let us look at Proposition \ref{N1inEmin} again. On one hand any
function $\varphi \in ({\cal N})_\mu ^1$ can be written in the form 
\begin{equation}
\label{phiPn}\varphi (z)=\sum_{n=0}^\infty \langle P_n^\mu (x),\varphi
^{(n)}\rangle \ ,\qquad \varphi ^{(n)}\in {\cal N}_{\Ckl }^{\hat \otimes n}\
, 
\end{equation}
on the other hand it is entire, i.e., it has the representation 
\begin{equation}
\label{phizn}\varphi (z)=\sum_{n=0}^\infty \langle z^{\otimes n},\tilde
\varphi ^{(n)}\rangle \ ,\qquad \tilde \varphi ^{(n)}\in {\cal N}_{\Ckl %
}^{\hat \otimes n}\ , 
\end{equation}
To proceed further we need the explicit correspondence $\left\{ \varphi
^{(n)},n\in \N \right\} \longleftrightarrow \left\{ \tilde \varphi
^{(n)},n\in \N \right\} $ which is given in the next lemma.

\begin{lemma}
\label{Reordering}{\bf (Reordering)} \smallskip
\\Equations (\ref{phiPn}) and (\ref{phizn}) hold iff%
$$
\tilde \varphi ^{(k)}=\sum_{n=0}^\infty \binom{n+k}k\left( P_n^\mu
(0),\varphi ^{(n+k)}\right) _{{\cal H}^{\hat \otimes n}} 
$$
or equivalently%
$$
\varphi ^{(k)}=\sum_{n=0}^\infty \binom{n+k}k\left( {\rm M}_n^\mu ,\tilde
\varphi ^{(n+k)}\right) _{{\cal H}^{\hat \otimes n}} 
$$
where $\left( P_n^\mu (0),\varphi ^{(n+k)}\right) _{{\cal H}^{\hat \otimes
n}}$ and $\left( {\rm M}_n^\mu ,\tilde \varphi ^{(n+k)}\right) _{{\cal H}%
^{\hat \otimes n}}$ denote contractions defined by (\ref{contraction}).
\end{lemma}

\noindent This is a consequence of (\ref{(P1)}) and (\ref{(P2)}). We omit
the simple proof.\bigskip\ \ 

Proposition \ref{N1inEmin} states%
$$
({\cal N})_\mu ^1\subseteq {\cal E}_{\min }^1({\cal N}^{\prime }) 
$$
as sets, where%
$$
{\cal E}_{\min }^1({\cal N}^{\prime })=\left\{ \varphi |_{{\cal N}^{\prime
}}\;\Big| \;\varphi \in {\cal E}_{\min }^1({\cal N}_{\Ckl }^{\prime
})\right\} \ . 
$$
Corollary \ref{phi(z)Betrag} then implies that the embedding is also
continuous. Now we are going to show that the converse also holds.

\begin{theorem}
\label{N1E1min}For all measures $\mu \in {\cal M}_a({\cal N}^{\prime })$ we
have the topological identity%
$$
({\cal N})_\mu ^1={\cal E}_{\min }^1({\cal N}^{\prime })\ . 
$$
\end{theorem}

\noindent To prove the missing topological inclusion it is convenient to use
the nuclear topology on ${\cal E}_{\min }^1({\cal N}_{\Ckl }^{\prime })$
(given by the norms $\lnorm \cdot \rnorm _{{p,q,1}}$) introduced in section 
\ref{Preliminaries}. Theorem \ref{Ekminprlim} ensures that this topology is
equivalent to the projective topology induced by the norms ${\rm n}_{p,l,k}$%
. Then the above theorem is an immediate consequence of the following norm
estimate.

\begin{proposition}
Let $p>p_\mu $ ($p_\mu $ as in Lemma \ref{equiLemma}) such that $\left\|
i_{p,p_\mu }\right\| _{HS}$ is finite and $q\in \N $ such that $2^{q/2}>K_p$
($K_p:=eC\left\| i_{p,p_\mu }\right\| _{HS}$ as in (\ref{MnmuNorm})). For
any $\varphi \in {\rm E}_{p,q}^1$ the restriction $\varphi |_{{\cal N}%
^{\prime }}$ is a function from $({\cal H}_p)_{q^{\prime },\mu }^1\ ,\
q^{\prime }<q$. Moreover the following estimate holds%
$$
||\varphi ||_{p,q^{\prime },\mu }\leq \lnorm  \varphi \rnorm
_{p,q,1}(1-2^{-q/2}K_p)^{-1}(1-2^{q^{\prime }-q})^{-1/2}\ . 
$$
\end{proposition}

\TeXButton{Proof}{\proof}Let $p,q\in \N $, $K_p$ be defined as above. A
function $\varphi \in {\rm E}_{p,q}^1$ has the representation (\ref{phizn}).
Using the Reordering lemma combined with (\ref{MnmuNorm}) and 
$$
\left| \tilde \varphi ^{(n)}\right| _p\leq \frac 1{n!}\,2^{-nq/2}\lnorm
\varphi \rnorm  _{p,q,1} 
$$
we obtain a representation of the form (\ref{phiPn}) where 
\begin{eqnarray*}
\left| \varphi ^{(n)}\right| _p & \leq & \sum_{k=0}^\infty \binom{n+k}k\left|%
{\rm M}_k^\mu \right| _{-p}\left| \tilde \varphi ^{(n+k)}\right| _p %
\\&\leq & \lnorm  \varphi \rnorm  _{p,q,1}\sum_{k=0}^\infty %
\binom{n+k}k\frac{k!}{(n+k)!}K_p^k\,2^{-(n+k)q/2} %
\\& \leq & \lnorm  \varphi \rnorm  _{p,q,1}%
\frac1{n!}2^{-nq/2}\sum_{k=0}^\infty (2^{-q/2}K_p)^k %
\\& \leq & \lnorm  \varphi \rnorm  _{p,q,1}%
\frac1{n!}2^{-nq/2}(1-2^{-q/2}K_p)^{-1} .%
\end{eqnarray*}
For $q^{\prime }<q$ this allows the following estimate%
\begin{eqnarray*}
||\varphi ||_{p,q^{\prime },\mu }^2%
&=&\sum_{n=0}^\infty (n!)^2\,2^{q^{\prime}n}\,|\varphi ^{(n)}|_p^2 %
\\& \leq & \lnorm  \varphi \rnorm ^2%
_{p,q,1}(1-2^{-q/2}K_p)^{-2}\sum_{k=0}^\infty 2^{n(q^{\prime }-q)}<\infty%
\end{eqnarray*}
This completes the proof.\TeXButton{End Proof}{\endproof}\bigskip\ 

Since we now have proved that the space of test functions $({\cal N})_\mu ^1$
is isomorphic to ${\cal E}_{\min }^1({\cal N}^{\prime })$ for all measures $%
\mu \in {\cal M}_a({\cal N}^{\prime })$, we will now drop the subscript $\mu 
$. The test function space $({\cal N})^1$ is the same for all measures $\mu
\in {\cal M}_a({\cal N}^{\prime })$.

\begin{corollary}
$({\cal N})^1$ is an algebra under pointwise multiplication.
\end{corollary}

\begin{corollary}
$({\cal N})^1$ admits `scaling' i.e., for $\lambda \in \C $ the scaling
operator $\sigma _\lambda :({\cal N})^1\rightarrow ({\cal N})^1$ defined by $%
\sigma _\lambda \varphi (x):=\varphi (\lambda x)$, $\varphi \in ({\cal N})^1$%
, $x\in {\cal N}^{\prime }$ is well--defined.
\end{corollary}

\begin{corollary}
For all $z\in {\cal N}_{\Ckl }^{\prime }$ the space $({\cal N})^1$ is
invariant under the shift operator $\tau _z:\varphi \mapsto \varphi (\cdot
+z)$.
\end{corollary}

\section{Distributions\label{Distributions}}

In this section we will introduce and study the space $({\cal N})_\mu ^{-1}$
of distributions corresponding to the space of test functions $({\cal N})^1$%
. Since ${\cal P}({\cal N}^{\prime })\subset ({\cal N})^1$ the space $({\cal %
N})_\mu ^{-1}$ can be viewed as a subspace of ${\cal P}_\mu ^{\prime }({\cal %
N}^{\prime })$%
$$
({\cal N})_\mu ^{-1}\subset {\cal P}_\mu ^{\prime }({\cal N}^{\prime }) 
$$
Let us now introduce the Hilbertian subspace $({\cal H}_{-p})_{-q,\mu }^{-1}$
of ${\cal P}_\mu ^{\prime }({\cal N}^{\prime })$ for which the norm 
$$
\left\| \Phi \right\| _{-p,-q,\mu }^2:=\sum_{n=0}^\infty 2^{-qn}\left| \Phi
^{(n)}\right| _{-p}^2\text{ } 
$$
is finite. Here we used the canonical representation%
$$
\Phi =\sum_{n=0}^\infty Q_n^\mu (\Phi ^{(n)})\in {\cal P}_\mu ^{\prime }(%
{\cal N}^{\prime })\text{ } 
$$
from Theorem \ref{PStrichRep}. The space $({\cal H}_{-p})_{-q,\mu }^{-1}$ is
the dual space of $({\cal H}_p)_q^1$ with respect to $L^2(\mu )$ (because of
the biorthogonality of $\p -$and $\Q -$systems). By general duality theory 
$$
({\cal N})_\mu ^{-1}:=\bigcup_{p,q\in \N }({\cal H}_{-p})_{-q,\mu }^{-1} 
$$
is the dual space of $({\cal N})^1$ with respect to $L^2(\mu )$. As we noted
in section \ref{Preliminaries} there exists a natural topology on co-nuclear
spaces (which coincides with the inductive limit topology). We will consider 
$({\cal N})_\mu ^{-1}$ as a topological vector space with this topology. So
we have the nuclear triple%
$$
({\cal N})^1\subset L^2(\mu )\subset ({\cal N})_\mu ^{-1}\ . 
$$
The action of $\Phi =\sum_{n=0}^\infty Q_n^\mu (\Phi ^{(n)})\in ({\cal N}%
)_\mu ^{-1}$ on a test function $\varphi =\sum_{n=0}^\infty \langle P_n^\mu
,\varphi ^{(n)}\rangle \in ({\cal N})^1$ is given by 
$$
\langle \!\langle \Phi ,\varphi \rangle \!\rangle _\mu =\sum_{n=0}^\infty
n!\langle \Phi ^{(n)},\varphi ^{(n)}\rangle \ . 
$$
\bigskip\ 

For a more detailed characterization of the singularity of distributions in $%
({\cal N})_\mu ^{-1}$ we will introduce some subspaces in this distribution
space. For $\beta \in [0,1]$ we define%
$$
({\cal H}_{-p})_{-q,\mu }^{-\beta }=\left\{ \Phi \in {\cal P}_\mu ^{\prime }(%
{\cal N}^{\prime })\;\bigg|\;\sum_{n=0}^\infty (n!)^{1-\beta }2^{-qn}\left|
\Phi ^{(n)}\right| _{-p}^2<\infty \text{ for }\Phi =\sum_{n=0}^\infty
Q_n^\mu (\Phi ^{(n)})\right\} 
$$
and 
$$
({\cal N})_\mu ^{-\beta }=\stackunder{p,q\in \N }{\bigcup }({\cal H}%
_{-p})_{-q,\mu }^{-\beta }\ , 
$$
It is clear that the singularity increases with increasing $\beta $:%
$$
({\cal N})^{-0}\subset ({\cal N})^{-\beta _1}\subset ({\cal N})^{-\beta
_2}\subset ({\cal N})^{-1} 
$$
if $\beta _1\leq \beta _2$.We will also consider $({\cal N})_\mu ^\beta $ as
equipped with the natural topology.\bigskip\ 

\example 
\newcounter{RadonNy} \setcounter{RadonNy}{\value{example}}{\it (Generalized
Radon--Nikodym derivative)} \smallskip 
\\We want to define a generalized function $\rho _\mu (z,\cdot )\in ({\cal N}%
)_\mu ^{-1}\ $, $z\in {\cal N}_{\Ckl }^{\prime }$ with the following property%
$$
\langle \!\langle \rho _\mu (z,\cdot ),\varphi \rangle \!\rangle _\mu =\int_{%
{\cal N}^{\prime }}\varphi (x-z)\;{\rm d}\mu (x)\ ,\qquad \varphi \in ({\cal %
N})^1\ . 
$$
That means we have to establish the continuity of $\rho _\mu (z,\cdot )$.
Let $z\in {\cal H}_{-p,\Ckl }$.\ If $p^{\prime }\geq p$ is sufficiently
large and $\varepsilon >0$ small enough, Corollary \ref{phi(z)Betrag}
applies i.e., $\exists q\in \N $ and $C>0$ such that%
\begin{eqnarray*}
\left| \int_{{\cal N}^{\prime }}\varphi (x-z){\rm d}\mu (x)\ \right| %
& \leq & C\left\| \varphi \right\| _{p^{\prime },q,\mu }%
\int_{{\cal N}^{\prime}}e^{\varepsilon |x-z|_{-p^{\prime }}}{\rm d}\mu (x) %
\\& \leq & C\left\| \varphi \right\| _{p^{\prime },q,\mu }e^{\varepsilon |z|_{-p^{\prime }}}\int_{{\cal N}^{\prime }}e^{\varepsilon |x|_{-p^{\prime }}}{\rm d}\mu (x) %
\end{eqnarray*}
If $\varepsilon $ is chosen sufficiently small the last integral exists.
Thus we have in fact $\rho (z,\cdot )\in ({\cal N})_\mu ^{-1}$. It is clear
that whenever the Radon--Nikodym derivative $\frac{{\rm d}\mu (x+\xi )}{{\rm %
d}\mu (x)}$ exists (e.g., $\xi \in {\cal N}$ in case $\mu $ is ${\cal N}$%
-quasi-invariant) it coincides with $\rho _\mu (\xi ,\cdot )$ defined above.
We will now show that in $({\cal N})_\mu ^{-1}$ we have the canonical
expansion%
$$
\rho _\mu (z,\cdot )=\sum_{n=0}^\infty \frac 1{n!}(-1)^nQ_n^\mu (z^{\otimes
n}). 
$$
It is easy to see that the r.h.s. defines an element in $({\cal N})_\mu
^{-1} $. Since both sides are in $({\cal N})_\mu ^{-1}$ it is sufficient to
compare their action on a total set from $({\cal N})^1$. For $\varphi
^{(n)}\in {\cal N}_{\Ckl }^{\hat \otimes n}$ we have 
\begin{eqnarray*}
\left\langle \!\!\left\langle \rho _\mu (z,\cdot ),\langle P_n^\mu ,\varphi%
^{(n)}\rangle \right\rangle \!\!\right\rangle _\mu %
&=& \int_{{\cal N}^{\prime }}\langle P_n^\mu (x-z),\varphi ^{(n)}\rangle%
 \;{\rm d}\mu (x) %
\\&=& \sum_{k=0}^\infty \binom nk(-1)^{n-k}\int_{{\cal N}^{\prime }}\langle %
P_k^\mu (x)\hat \otimes z^{\otimes n-k},\varphi ^{(n)}\rangle \;{\rm d}\mu %
(x) \\&=& (-1)^n\langle z^{\otimes n},\varphi ^{(n)}\rangle %
\\&=& \left\langle \!\!\left\langle \sum_{k=0}^\infty \frac 1{k!}(-1)^kQ_k^\mu %
(z^{\otimes k}),\langle P_n^\mu ,\varphi ^{(n)}\rangle \right\rangle%
\!\!\right\rangle _\mu \ ,%
\end{eqnarray*}
where we have used (\ref{(P3)}), (\ref{(P4)}) and the biorthogonality of $%
\p 
$- and $\Q $-systems. This had to be shown. In other words, we have proven
that $\rho _\mu (-z,\cdot )$ is the generating function of the $\Q $%
-functions 
\begin{equation}
\label{rhomyQn}\rho _\mu (-z,\cdot )=\sum_{n=0}^\infty \frac 1{n!}Q_n^\mu
(z^{\otimes n})\ . 
\end{equation}
Let use finally remark that the above expansion allows for more detailed
estimates. It is easy to see that $\rho _\mu \in ({\cal N})_\mu ^{-0}$.%
\bigskip\ 

\example {\it (Delta distribution)} \smallskip \\For $z\in {\cal N}_{\Ckl %
}^{\prime }$ we define a distribution by the following $\Q $-decomposition:%
$$
\delta _z=\sum_{n=0}^\infty \frac 1{n!}Q_n^\mu (P_n^\mu (z)) 
$$
If $p\in \N $ is large enough and $\varepsilon >0$ sufficiently small there
exists $C_{p,\varepsilon }>0$ according to (\ref{(P5)}) such that 
\begin{eqnarray*}
\left\| \delta _z\right\| _{-p,-q,\mu }^2 &=&\sum_{n=0}^\infty %
(n!)^{-2}2^{-nq}\left| P_n^\mu (z)\right| _{-p}^2 %
\\& \leq & C_{p,\varepsilon }^2\,e^{2\varepsilon |z|_{-p}}\sum_{n=0}^\infty%
2^{-nq}\varepsilon ^{-2n}\ ,\qquad z\in {\cal H}_{-p,\Ckl }\;, %
\end{eqnarray*}
which is finite for sufficiently large $q\in \N $. Thus $\delta _z\in ({\cal %
N})_\mu ^{-1}$.

For $\varphi =\sum_{n=0}^\infty \langle P_n^\mu ,\varphi ^{(n)}\rangle \in (%
{\cal N})^1$ the action of $\delta _z$ is given by 
$$
\langle \!\langle \delta _z,\varphi \rangle \!\rangle _\mu
=\sum_{n=0}^\infty \langle P_n^\mu (z),\varphi ^{(n)}\rangle =\varphi (z) 
$$
because of (\ref{QnPnPair}). This means that $\delta _z$ (in particular for $%
z$ real) plays the role of a ``$\delta $-function" (evaluation map) in the
calculus we discuss.

\section{Integral transformations}

\begin{sloppypar}
We will first introduce the Laplace transform of a function $\varphi \in
L^2(\mu )$. The global assumption $\mu \in {\cal M}_a({\cal N}^{\prime })$
guarantees the existence of $p_\mu ^{\prime }\in \N \ $, $\varepsilon _\mu
>0 $ such that $\int_{{\cal N}^{\prime }}\exp (\varepsilon _\mu |x|_{-p_\mu
^{\prime }})\,{\rm d}\mu (x)<\infty $ by Lemma \ref{equiLemma}. Thus $\exp
(\langle x,\theta \rangle )\in L^2(\mu )$ if $2|\theta |_{p_\mu ^{\prime
}}\leq \varepsilon _\mu \ ,\theta \in {\cal H}_{p_\mu ^{\prime },\Ckl }$.
Then by Cauchy--Schwarz inequality the Laplace transform defined by 
$$
L_\mu \varphi (\theta ):=\int_{{\cal N}^{\prime }}\varphi (x)\exp \langle
x,\theta \rangle \,{\rm d}\mu (x) 
$$
is well defined for $\varphi \in L^2(\mu )\ ,\theta \in {\cal H}_{p_\mu
^{\prime },\Ckl }$ with $2|\theta |_{p_\mu ^{\prime }}\leq \varepsilon _\mu $%
. Now we are interested to extend this integral transform from $L^2(\mu )$
to the space of distributions $({\cal N})_\mu ^{-1}$.
\end{sloppypar}

Since our construction of test function and distribution spaces is closely
related to $\p $- and $\Q $-systems it is useful to introduce the so called $%
S_\mu $-transform 
$$
S_\mu \varphi (\theta ):=\frac{L_\mu \varphi (\theta )}{l_\mu (\theta )}\ . 
$$
Since $e_\mu (\theta ;x)=e^{\langle x,\theta \rangle }/l_\mu (\theta )$ we
may also write%
$$
S_\mu \varphi (\theta )=\int_{{\cal N}^{\prime }}\varphi (x)\,e_\mu (\theta
;x)\,{\rm d}\mu (x)\ . 
$$
The $\mu $-exponential $e_\mu (\theta ,\cdot )$ is not a test function in $(%
{\cal N})^1$, see Example \arabic{myexponent} . So the definition of the $%
S_\mu $-transform of a distribution $\Phi \in ({\cal N})_\mu ^{-1}$ must be
more careful. Every such $\Phi $ is of finite order i.e., $\exists p,q\in 
\N 
$ such that $\Phi \in ({\cal H}_{-p})_{-q,\mu .}^{-1}$ As shown in Example 
\arabic{myexponent} $e_\mu (\theta ,\cdot )$ is in the corresponding dual
space $({\cal H}_p)_{q,\mu }^1$ if $\theta \in {\cal H}_{p,\Ckl }$ is such
that $2^q|\theta |_p^2<1$. Then we can define a consistent extension of $%
S_\mu $-transform.%
$$
S_\mu \Phi (\theta ):=\langle \!\langle \Phi ,\;e_\mu (\theta ,\cdot
)\rangle \!\rangle _\mu 
$$
if $\theta $ is chosen in the above way. The biorthogonality of $\p $- and $%
\Q $-system implies 
$$
S_\mu \Phi (\theta )=\sum_{n=0}^\infty \langle \Phi ^{(n)},\theta ^{\otimes
n}\rangle \ . 
$$
It is easy to see that the series converges uniformly and absolutely on any
closed ball $\left\{ \left. \theta \in {\cal H}_{p,\Ckl }\right| \;|\theta
|_p^2\leq r,\ r<2^{-q}\right\} $, see the proof of Theorem \ref{CharTh}.
Thus $S_\mu \Phi $ is holomorphic a neighborhood of zero, i.e., $S_\mu \Phi
\in {\rm Hol}_0({\cal N}_{\Ckl })$. In the next section we will discuss this
relation to the theory of holomorphic functions in more detail.

The third integral transform we are going to introduce is more appropriate
for the test function space $({\cal N})^1$. We introduce the convolution of
a function $\varphi \in ({\cal N})^1$ with the measure $\mu $ by%
$$
C_\mu \varphi (y):=\int_{{\cal N}^{\prime }}\varphi (x+y)\,{\rm d}\mu (x) ,
\quad y\in {\cal N}^{\prime } .%
$$
From Example \arabic{RadonNy} the existence of a generalized Radon--Nikodym
derivative $\rho _\mu (z,\cdot )$, $z\in {\cal N}_{\Ckl }^{\prime }$ in $(%
{\cal N})_\mu ^{-1}$ is guaranteed. So for any $\varphi \in ({\cal N})^1$, $%
z\in {\cal N}_{\Ckl }^{\prime }$ the convolution has the representation 
$$
C_\mu \varphi (z)=\langle \!\langle \rho _\mu (-z,\cdot ),\;\varphi \rangle
\!\rangle _\mu \;. 
$$
If $\varphi \in ({\cal N})^1$ has the canonical representation%
$$
\varphi =\sum_{n=0}^\infty \langle P_n^\mu ,\;\varphi ^{(n)}\rangle 
$$
we have by equation (\ref{rhomyQn}) 
$$
C_\mu \varphi (z)=\sum_{n=0}^\infty \langle z^{\otimes n},\varphi
^{(n)}\rangle \ . 
$$

In Gaussian Analysis $C_\mu $- and $S_\mu$-transform coincide. It is a
typical non-Gaussian effect that these two transformations differ from each
other.

\section{Characterization theorems \label{Characterization}}

Gaussian Analysis has shown that for applications it is very useful to
characterize test and distribution spaces by the integral transforms
introduced in the previous section. In the non-Gaussian setting first
results in this direction have been obtained by \cite{AKS93,ADKS94}.%
\bigskip 

We will start to characterize the space $({\cal N})^1$ in terms of the
convolution $C_\mu $.

\begin{theorem}
\label{CmuChar} The convolution $C_\mu $ is a topological isomorphism from $(%
{\cal N})^1$ on ${\cal E}_{\min }^1({\cal N}_{\Ckl }^{\prime })$.
\end{theorem}

\TeXButton{Remark }{\remark } Since we have identified $({\cal N})^1$ and $%
{\cal E}_{\min }^1({\cal N}^{\prime })$ by Theorem \ref{N1E1min} the above
assertion can be restated as follows. We have 
$$
C_\mu :{\cal E}_{\min }^1({\cal N}^{\prime })\rightarrow {\cal E}_{\min }^1(%
{\cal N}_{\Ckl }^{\prime }) 
$$
as a topological isomorphism.

\TeXButton{Proof}{\proof}The proof has been well prepared by Theorem \ref
{Ekminprlim}, because the nuclear topology on ${\cal E}_{\min }^1({\cal N}_{%
\Ckl }^{\prime })$ is the most natural one from the point of view of the
above theorem. Let $\varphi \in ({\cal N})^1$ with the representation%
$$
\varphi =\sum_{n=0}^\infty \langle P_n^\mu ,\varphi ^{(n)}\rangle \ . 
$$
From the previous section it follows%
$$
C_\mu \varphi (z)=\sum_{n=0}^\infty \langle z^{\otimes n},\varphi
^{(n)}\rangle \ 
$$
It is obvious from (\ref{3StrichNorm}) that%
$$
\lnorm  C_\mu \varphi \rnorm  _{p,q,1}=\left\| \varphi \right\| _{p,q,\mu }\ 
$$
for all $p,q\in \N _0$, which proves the continuity of 
$$
C_\mu :({\cal N})^1\rightarrow {\cal E}_{\min }^1({\cal N}_{\Ckl }^{\prime
})\ . 
$$

Conversely let $F\in {\cal E}_{\min }^1({\cal N}_{\Ckl }^{\prime })$. Then
Theorem \ref{Ekminprlim} ensures the existence of a sequence of generalized
kernels $\left\{ \varphi ^{(n)}\in {\cal N}_{\Ckl }^{\prime }\;|\;n\in \N %
_0\right\} $ such that 
$$
F(z)=\sum_{n=0}^\infty \langle z^{\otimes n},\varphi ^{(n)}\rangle \ . 
$$
Moreover for all $p,q\in \N _0$%
$$
\lnorm  F\rnorm  _{p,q,1}^2=\sum_{n=0}^\infty (n!)^2\,2^{nq}\left| \varphi
^{(n)}\right| _p^2 
$$
is finite. Choosing%
$$
\varphi =\sum_{n=0}^\infty \langle P_n^\mu ,\varphi ^{(n)}\rangle 
$$
we have $\left\| \varphi \right\| _{p,q,\mu }=\lnorm  F\rnorm
_{p,q,1}$. Thus $\varphi \in ({\cal N})^1$. Since $C_\mu \varphi =F$ we have
shown the existence and continuity of the inverse of $C_\mu $.%
\TeXButton{End Proof}{\endproof}\bigskip\ 

To illustrate the above theorem in terms of the natural topology on ${\cal E}%
_{\min }^1({\cal N}_{\Ckl }^{\prime })$ we will reformulate the above
theorem and add some useful estimates which relate growth in ${\cal E}_{\min
}^1({\cal N}_{\Ckl }^{\prime })$ to norms on $({\cal N})^1$.

\begin{corollary}
\hfill \\1) Let $\varphi \in ({\cal N})^1$ then for all $p,l\in \N _0$ and $%
z\in {\cal H}_{-p,\Ckl }$ the following estimate holds%
$$
\left| C_\mu \varphi (z)\right| \leq \left\| \varphi \right\| _{p,2l,\mu
}\exp (2^{-l}|z|_{-p}) 
$$
i.e., C$_\mu \varphi \in {\cal E}_{\min }^1({\cal N}_{\Ckl }^{\prime })$.%
\medskip\ \\2) Let $F\in {\cal E}_{\min }^1({\cal N}_{\Ckl }^{\prime })$.
Then there exists $\varphi \in ({\cal N})^1$ with $C_\mu \varphi =F$. The
estimate 
$$
\left| F(z)\right| \leq C\exp (2^{-l}|z|_{-p}) 
$$
for $C>0,\ p,q\in \N _0$ implies%
$$
\left\| \varphi \right\| _{p^{\prime },q,\mu }\leq C\left(
1-2^{q-2l}e^2\left\| i_{p^{\prime },p}\right\| _{HS}^2\right) ^{-1/2} 
$$
if the embedding $i_{p^{\prime },p}:{\cal H}_{p^{\prime }}\hookrightarrow 
{\cal H}_p$ is Hilbert-Schmidt and $2^{l-q/2}>e\left\| i_{p^{\prime
},p}\right\| _{HS}$.
\end{corollary}

\TeXButton{Proof}{\proof}The first statement follows from 
$$
\left| C_\mu \varphi (z)\right| \leq {\rm n}_{p,l,1}(C_\mu \varphi )\cdot
\exp (2^{-l}|z|_{-p}) 
$$
which follows from the definition of {\rm n}$_{p,l,1}$ and estimate (\ref
{nplk3StrichNorm}). The second statement is an immediate consequence of
Lemma \ref{3Strichnplk}. \TeXButton{End Proof}{\endproof}\bigskip\ 

The next theorem characterizes distributions from $({\cal N})_\mu ^{-1}$ in
terms of $S_\mu$-transform.

\begin{theorem}
\label{CharTh}The $S_\mu $-transform is a topological isomorphism from $(%
{\cal N})_\mu ^{-1}$ on ${\rm Hol}_0({\cal N}_{\Ckl })$.
\end{theorem}

\TeXButton{Remark }{\remark } The above theorem is closely related to the
second part of Theorem \ref{indlimEkmax}. Since we left the proof open we
will give a detailed proof here.

\TeXButton{Proof}{\proof}Let $\Phi \in ({\cal N})_\mu ^{-1}$ . Then there
exists $p,q\in \N $ such that 
$$
\left\| \Phi \right\| _{-p,-q,\mu }^2=\sum_{n=0}^\infty 2^{-nq}|\Phi
^{(n)}|_{-p}^2 
$$
is finite. From the previous section we have 
\begin{equation}
\label{SPhitheta}S_\mu \Phi (\theta )=\sum_{n=0}^\infty \langle \Phi
^{(n)},\theta ^{\otimes n}\rangle \;. 
\end{equation}
For $\theta \in {\cal N}_{\Ckl }$ such that $2^q|\theta |_p^2<1$ we have by
definition (Formula (\ref{3StrichNorm}))%
$$
\lnorm  S_\mu \Phi \rnorm  _{-p,-q,-1}=\left\| \Phi \right\| _{-p,-q,\mu \
}. 
$$
By Cauchy--Schwarz inequality 
\begin{eqnarray*}
\left| S_\mu \Phi (\theta )\right| & \leq & \sum_{n=0}^\infty |\Phi %
^{(n)}|_{-p}|\theta |_p^n %
\\&\leq & \left( \sum_{n=0}^\infty 2^{-nq}|\Phi ^{(n)}|_{-p}^2\right)%
^{1/2}\left( \sum_{n=0}^\infty 2^{nq}|\theta |_p^{2n}\right) ^{1/2} %
\\&=& \left\| \Phi \right\| _{-p,-q,\mu }\left( 1-2^q|\theta |_p^2\right)%
^{-1/2}\ . 
\end{eqnarray*}
Thus the series (\ref{SPhitheta}) converges uniformly on any closed ball $%
\left\{ \left. \theta \in {\cal H}_{p,\Ckl }\right| \;|\theta |_p^2\leq r,\
r<2^{-q}\right\} $. Hence $S_\mu \Phi \in {\rm Hol}_0({\cal N}_{\Ckl })$ and 
$$
{\rm n}_{p,l,\infty }(S_\mu \Phi )\leq \left\| \Phi \right\| _{-p,-q,\mu
}(1-2^{q-2l})^{-1/2} 
$$
if $2l>q$. This proves that $S_\mu $ is a continuous mapping from $({\cal N}%
)_\mu ^{-1}$ to ${\rm Hol}_0({\cal N}_{\Ckl })$. In the language of section 
\ref{Holomorphy} this reads%
$$
\stackunder{p,q\in \N }{\rm ind\ lim}\,{\rm E}_{-p,-q}^{-1}\subset {\rm Hol}%
_0({\cal N}_{\Ckl }) 
$$
topologically.\medskip\ 

Conversely, let $F\in {\rm Hol}_0({\cal N}_{\Ckl })$ be given, i.e., there
exist $p,l\in \N $ such that {\rm n}$_{p,l,\infty }(F)<\infty $. The first
step is to show that there exists $p^{\prime },q\in \N $ such that 
$$
\lnorm  F\rnorm  _{-p^{\prime },-q,-1}<{\rm n}_{p,l,\infty }(F)\cdot C\ , 
$$
for sufficiently large $C>0$. This implies immediately 
$$
{\rm Hol}_0({\cal N}_{\Ckl })\subset \ \stackunder{p,q\in \N }{\rm ind\ lim}%
\,{\rm E}_{-p,-q}^{-1} 
$$
topologically, which is the missing part in the proof of the second
statement in Theorem \ref{indlimEkmax}.

By assumption the Taylor expansion%
$$
F(\theta )=\sum_{n=0}^\infty \frac 1{n!}\widehat{{\rm d}^nF(0)}(\theta ) 
$$
converges uniformly on any closed ball $\left\{ \left. \theta \in {\cal H}%
_{p,\Ckl }\right| \;|\theta |_p^2\leq r,\ r<2^{-l}\right\} $ and 
$$
\left| F(\theta )\right| \leq {\rm n}_{p,l,\infty }(F)\ . 
$$
Proceeding analogously to Lemma \ref{nplk3Strich}, an application of
Cauchy's inequality gives%
\begin{eqnarray*}
\frac 1{n!}\widehat{{\rm d}^nF(0)}(\theta ) & \leq & 2^l|\theta |_p^n\sup %
_{|\theta |_p\leq 2^{-l}}|F(\theta )| %
\\& \leq & {\rm n}_{p,l,\infty }(F)\ \cdot 2^{nl}\cdot |\theta |_p^n 
\end{eqnarray*}
The polarization identity gives 
$$
\left| \frac 1{n!}{\rm d}^nF(0)(\theta _1,\ldots ,\theta _n)\right| \leq 
{\rm n}_{p,l,\infty }(F)\ \cdot e^n\cdot 2^{nl}\prod_{j=1}^n|\theta _j|_p 
$$
Then by kernel theorem (Theorem \ref{KernelTh}) there exist kernels $\Phi
^{(n)}\in {\cal H}_{-p^{\prime },\Ckl }^{\hat \otimes n}$ for $p^{\prime }>p$
with $\left\| i_{p^{\prime },p}\right\| _{HS}<\infty $ such that%
$$
F(\theta )=\sum_{n=0}^\infty \langle \Phi ^{(n)},\theta ^{\otimes n}\rangle
\ . 
$$
Moreover we have the following norm estimate 
$$
\left| \Phi ^{(n)}\right| _{-p^{\prime }}\leq {\rm n}_{p,l,\infty }(F)\
\left( 2^le\left\| i_{p^{\prime },p}\right\| _{HS}\right) ^n 
$$
Thus%
\begin{eqnarray*}
\lnorm  F\rnorm  _{-p^{\prime },-q,-1}^2%
&=& \sum_{n=0}^\infty%
2^{-nq}\left| \Phi ^{(n)}\right| _{-p^{\prime }}^2 %
\\& \leq & {\rm n}_{p,l,\infty }^2(F)\sum_{n=0}^\infty %
\left( 2^{2l-q}e^2\left\|%
i_{p^{\prime },p}\right\| _{HS}^2\right) ^n %
\\&=& {\rm n}_{p,l,\infty }^2(F)\left( 1-2^{2l-q}e^2%
\left\| i_{p^{\prime },p}\right\| _{HS}^2\right) ^{-1}\ 
\end{eqnarray*}
if $q\in \N $ is such that $\rho :=2^{2l-q}e^2\left\| i_{p^{\prime
},p}\right\| _{HS}^2\ <1$. So we have in fact 
$$
\lnorm  F\rnorm  _{-p^{\prime },-q,-1}\leq {\rm n}_{p,l,\infty }(F)(1-\rho
)^{-1/2}. 
$$
Now the rest is simple. Define $\Phi \in ({\cal N})_\mu ^{-1}$ by 
$$
\Phi =\sum_{n=0}^\infty Q_n^\mu (\Phi ^{(n)}) 
$$
then $S_\mu \Phi =F$ and 
$$
\left\| \Phi \right\| _{-p^{\prime },-q,\mu }=\lnorm  F\rnorm
_{-p^{\prime },-q,-1} 
$$
This proves the existence of a continuous inverse of the $S_\mu $%
--transform. Uniqueness of $\Phi $ follows from the fact that $\mu $%
-exponentials are total in any $({\cal H}_p)_q^1$. \TeXButton{End Proof}
{\endproof}\bigskip\ 

We can extract some useful estimates from the above proof which describe the
degree of singularity of a distribution.

\begin{corollary}
Let $F\in {\rm Hol}_0({\cal N}_{\Ckl })$ be holomorphic for all $\theta \in 
{\cal N}_{\Ckl }$ with $|\theta |_p\leq 2^{-l}$. If $p^{\prime }>p$ with $%
\left\| i_{p^{\prime },p}\right\| _{HS}<\infty $ and $q\in \N $ is such that 
$\rho :=2^{2l-q}e^2\left\| i_{p^{\prime },p}\right\| _{HS}^2<1$. Then $\Phi
\in ({\cal H}_{-p^{\prime }})_{-q}^{-1}$ and 
$$
\left\| \Phi \right\| _{-p^{\prime },-q,\mu }\leq {\rm n}_{p,l,\infty
}(F)\cdot (1-\rho )^{-1/2}. 
$$
\end{corollary}

For a more detailed discussion of the degree of singularity the spaces $(%
{\cal N})^{-\beta },\ \beta \in [0,1)$ are useful. In the following theorem
we will characterize these spaces by means of $S_\mu$-transform.

\begin{theorem}
\sloppy The $S_\mu $-transform is a topological isomorphism from $({\cal N}%
)_\mu ^{-\beta }$, $\beta \in [0,1)$ on ${\cal E}_{\max }^{2/(1-\beta )}(%
{\cal N}_{\Ckl })$.
\end{theorem}

\fussy \TeXButton{Remark }{\remark } The proof will also complete the proof
of Theorem \ref{indlimEkmax}.

\TeXButton{Proof}{\proof}Let $\Phi \in ({\cal H}_{-p})_{-q,\mu }^{-\beta }$
with the canonical representation $\Phi =\sum_{n=0}^\infty Q_n^\mu (\Phi
^{(n)})$ be given. The $S_\mu $-transform of $\Phi $ is given by 
$$
S_\mu \Phi (\theta ) =\sum_{n=0}^\infty \langle \Phi
^{(n)},\theta ^{\otimes n}\rangle . 
$$
Hence 
$$
\lnorm  S_\mu \Phi \rnorm  _{-p,-q,-\beta }^2=\sum_{n=0}^\infty
(n!)^{1-\beta }\,2^{-nq}|\Phi ^{(n)}|_{-p}^2 
$$
is finite. We will show that there exist $l\in \N $ and $C<0$ such that%
$$
{\rm n}_{-p,-l,2/(1-\beta )}(S_\mu \Phi )\leq C\,\lnorm  S_\mu \Phi \rnorm  %
_{-p,-q,-\beta }\ . 
$$

We can estimate as follows%
\begin{eqnarray*}
|S_\mu \Phi (\theta )| & \leq & \sum_{n=0}^\infty \left| \Phi ^{(n)}\right| %
_{-p}\left| \theta \right| _p^n %
\\& \leq & \left( \sum_{n=0}^\infty (n!)^{1-\beta }2^{-nq}|\Phi %
^{(n)}|_{-p}^2\right) ^{1/2}\left( \sum_{n=0}^\infty %
\frac 1{(n!)^{1-\beta }}2^{nq}\left| \theta \right| _p^{2n}\right) ^{1/2} %
\\&=& \lnorm  S_\mu \Phi \rnorm  _{-p,-q,-\beta }\left( %
\sum_{n=0}^\infty \rho ^{n\beta }\cdot \frac 1{(n!)^{1-\beta }}2^{nq}\rho %
^{-n\beta }\left| \theta \right| _p^{2n}\cdot \right) ^{1/2}, 
\end{eqnarray*}
where we have introduced a parameter $\rho \in (0,1)$. An application of
H\"older's inequality for the conjugate indices $\frac 1\beta $ and $\frac
1{1-\beta }$ gives%
\begin{eqnarray*}
\left| S_\mu \Phi (\theta )\right| & \leq & \lnorm  S_\mu \Phi \rnorm%
_{-p,-q,-\beta }\left( \sum_{n=0}^\infty \rho ^n\right) %
^{\beta /2}\cdot \left( \sum_{n=0}^\infty \frac 1{n!}\left( 2^q\rho %
^{-\beta }|\theta |_p^2\right) ^{\frac n{1-\beta }}\right) %
^{\frac{1-\beta }2} %
\\&=& \lnorm  S_\mu \Phi \rnorm  %
_{-p,-q,-\beta }\left( 1-\rho %
\right) ^{-\beta /2}\exp \left( \tfrac{1-\beta }2 \,2^{\frac q{1-\beta }}\,%
\rho ^{-\frac \beta {1-\beta }}\, |\theta |_p^{\frac 2{1-\beta }}\right) %
\end{eqnarray*}
If $l\in \N $ is such that%
$$
2^{l-\frac q{1-\beta }}>\tfrac{1-\beta }2\rho ^{-\frac \beta {1-\beta }} 
$$
we have%
\begin{eqnarray*}
{\rm n}_{-p,-l,2/(1-\beta )}(S_\mu \Phi )%
&=& \sup _{\theta \in {\cal H}_{p,\Ckkl }}\left| S_\mu \Phi (\theta )%
\right| \, \exp \left( -2^l|\theta %
|_p^{2/(1-\beta )}\right) %
\\& \leq & \left( 1-\rho \right) ^{-\beta /2}\lnorm  S_\mu \Phi \rnorm %
_{-p,-q,-\beta }%
\end{eqnarray*}
This shows that $S_\mu $ is continuous from $({\cal N})_\mu ^{-\beta }$ to $%
{\cal E}_{\max }^{2/(1-\beta )}({\cal N}_{\Ckl }).$ Or in the language of
Theorem \ref{indlimEkmax}%
$$
\stackunder{p,q\in \N }{\rm ind\ lim}\,{\rm E}_{-p,-q}^{-\beta }\subset 
{\cal E}_{\max }^{2/(1-\beta )}({\cal N}_{\Ckl }) 
$$
topologically.\medskip\ 

The proof of the inverse direction is closely related to the proof of Lemma 
\ref{3Strichnplk}. So we will be more sketchy in the following.

\noindent Let $F\in {\cal E}_{\max }^k({\cal N}_{\Ckl }),\ k=\frac 2{1-\beta
}$. Hence there exist $p,l\in \N _0$ such that%
$$
\left| F(\theta )\right| \leq {\rm n}_{-p,-l,k}(F)\exp (2^l|\theta |_p^k)\
,\qquad \theta \in {\cal N}_{\Ckl } 
$$
From this we have completely analogous to the proof of Lemma \ref
{3Strichnplk} by Cauchy inequality and kernel theorem the representation%
$$
F(\theta )=\sum_{n=0}^\infty \langle \Phi ^{(n)},\theta ^{\otimes n}\rangle 
$$
and the bound%
$$
\left| \Phi ^{(n)}\right| _{-p^{\prime }}\leq {\rm n}_{-p,-l,k}(F)%
\;(n!)^{-1/k}\left\{ (k2^l)^{1/k}e\left\| i_{p^{\prime },p}\right\|
_{HS}\right\} ^n\ , 
$$
where $p^{\prime }>p$ is such that $i_{p^{\prime },p}:{\cal H}_{p^{\prime
}}\hookrightarrow {\cal H}_p$ is Hilbert--Schmidt. Using this we have 
\begin{eqnarray*}
\lnorm  F\rnorm  _{-p^{\prime },-q,-\beta}^2 %
&=&  \sum_{n=0}^\infty (n!)^{1-\beta }2^{-qn}\left| \Phi ^{(n)}\right| %
_{-p^{\prime }}^2 %
\\& \leq & {\rm n}_{-p,-l,k}^2(F)\sum_{n=0}^\infty (n!)^{1-\beta -2/k}2^{-qn}%
\left\{ (k2^l)^{1/k}e\left\| i_{p^{\prime },p}\right\| %
_{HS}\right\} ^{2n} %
\\& \leq & {\rm n}_{-p,-l,k}^2(F)\sum_{n=0}^\infty \rho ^n 
\end{eqnarray*}
where we have set $\rho :=2^{-q+2l/k}k^{2/k}e^2\left\| i_{p^{\prime
},p}\right\| _{HS}^2$ . If $q\in \N $ is chosen large enough such that $\rho
<1$ the sum on the right hand side is convergent and we have 
\begin{equation}
\label{F3Strich}\lnorm  F\rnorm  _{-p^{\prime },-q,-\beta }\leq {\rm n}%
_{-p,-l,2/(1-\beta )}(F)\cdot (1-\rho )^{-1/2}\ . 
\end{equation}
That means%
$$
{\cal E}_{\max }^{2/(1-\beta )}({\cal N}_{\Ckl })\subset \ \stackunder{%
p,q\in \N _0}{\rm ind\ lim}\,{\rm E}_{-p,-q}^{-\beta } 
$$
topologically.

If we set 
$$
\Phi :=\sum_{n=0}^\infty Q_n^\mu (\Phi ^{(n)}) 
$$
then $S_\mu\Phi =F$ and $\Phi \in ({\cal H}_{-p^{\prime }})_{-q}^{-\beta }$
since%
$$
\sum_{n=0}^\infty (n!)^{1-\beta }2^{-qn}|\Phi ^{(n)}|_{-p^{\prime }}^2 
$$
is finite. Hence%
$$
S_\mu:({\cal N})_\mu ^{-\beta }\rightarrow {\cal E}_{\max }^{2/(1-\beta )}(%
{\cal N}_{\Ckl }) 
$$
is one to one. The continuity of the inverse mapping follows from the norm
estimate (\ref{F3Strich}).\TeXButton{End Proof}{\endproof}

%% file: wick.tex
\LaTeXparent{nga4.tex}

\section{The Wick product \label{Wick}}

In Gaussian Analysis it has been shown that $({\cal N})_{\gamma _{{\cal H}%
}}^{-1}$ (and other distribution spaces) is closed under so called Wick
multiplication (see \cite{KLS94} and \cite{BeS95,Ok94,Va95} for
applications). This concept has a natural generalization to the present
setting.

\begin{definition}
\ Let $\Phi ,\Psi \in $ $({\cal N})_\mu ^{-1}$. Then we define the Wick
product $\Phi \diamond \Psi $by 
$$
S_\mu (\Phi \diamond \Psi )=S_\mu \Phi \cdot S_\mu \Psi \ . 
$$
\end{definition}

This is well defined because ${\rm Hol}_0({\cal N}_{\Ckl })$ is an algebra
and thus by the characterization Theorem \ref{CharTh} there exists an
element $\Phi \diamond \Psi \in ({\cal N})_\mu ^{-1}$ such that $S_\mu (\Phi
\diamond \Psi )=S_\mu \Phi \cdot S_\mu \Psi $.

By definition we have 
$$
Q_n^\mu (\Phi ^{(n)})\diamond Q_m^\mu (\Psi ^{(m)})=Q_{n+m}^\mu (\Phi
^{(n)}\hat \otimes \Psi ^{(m)})\text{ ,} 
$$
$\Phi ^{(n)}\in ({\cal N}_{\Ckl }^{\hat \otimes n})^{\prime }$ and $\Psi
^{(m)}\in ({\cal N}_{\Ckl }^{\hat \otimes m})^{\prime }$. So in terms of $%
\Q 
$--decompositions $\Phi =\sum_{n=0}^\infty Q_n^\mu (\Phi ^{(n)})$ and $\Psi
=\sum_{n=0}^\infty Q_n^\mu (\Psi ^{(n)})$ the Wick product is given by 
$$
\Phi \diamond \Psi =\sum_{n=0}^\infty Q_n^\mu (\Xi ^{(n)}) 
$$
where%
$$
\Xi ^{(n)}=\sum_{k=0}^n\Phi ^{(k)}\hat \otimes \Psi ^{(n-k)} 
$$
This allows for concrete norm estimates.

\begin{proposition}
The Wick product is continuous on $({\cal N})_\mu ^{-1}$. In particular the
following estimate holds for $\Phi \in ({\cal H}_{-p_1})_{-q_1,\mu }^{-1}\
,\ \Psi \in ({\cal H}_{-q_2})_{-q_2}^{-1}$ and $p=\max (p_1,p_2),\
q=q_1+q_2+1$ 
$$
\left\| \Phi \diamond \Psi \right\| _{-p,-q,\mu }=\left\| \Phi \right\|
_{-p_1,-q_1,\mu }\left\| \Psi \right\| _{-p_2,-q_2,\mu }\text{ .} 
$$
\end{proposition}

\TeXButton{Proof}{\proof}We can estimate as follows%
\begin{eqnarray*}
\left\| \Phi \diamond \Psi \right\| _{-p,-q,\mu }^2 %
&=& \sum_{n=0}^\infty 2^{-nq}\left| \Xi ^{(n)}\right| _{-p}^2 %
\\&=& \sum_{n=0}^\infty 2^{-nq}\left( \sum_{k=0}^n\left| \Phi ^{(k)}\right| %
_{-p}\left| \Psi ^{(n-k)}\right| _{-p}\right) ^2 %
\\& \leq & \sum_{n=0}^\infty 2^{-nq}\,(n+1)\sum_{k=0}^n\left| \Phi ^{(k)} %
\right| _{-p}^2\left| \Psi ^{(n-k)}\right| _{-p}^2 %
\\& \leq & \sum_{n=0}^\infty \sum_{k=0}^n2^{-nq_1}\left| \Phi ^{(n)}\right| %
_{-p}^22^{-nq_2}\left| \Psi ^{(n-k)}\right| _{-p}^2 %
\\& \leq & \left( \sum_{n=0}^\infty 2^{-nq_1}\left| \Phi ^{(k)}\right| %
_{-p_1}^2\right) \left( \sum_{n=0}^\infty 2^{-nq_2}\left| \Psi ^{(n)}\right| %
_{-p_2}^2\right) %
\\&=& \left\| \Phi \right\| _{-p_1,-q_1,\mu }^2\left\| \Psi \right\|%
_{-p_2,-q_2,\mu }^2\text{ .} 
\end{eqnarray*}
\TeXButton{End Proof}{\endproof}\bigskip\ 

Similar to the Gaussian case the special properties of the space $({\cal N}%
)_\mu ^{-1}$ allow the definition of {\it Wick analytic functions }under
very general assumptions. This has proven to be of some relevance to solve
equations e.g., of the type $\Phi \diamond X=\Psi $ for $X\in ({\cal N})_\mu
^{-1}$ . See \cite{KLS94} for the Gaussian case.

\begin{theorem}
Let $F:\C \rightarrow \C $ be analytic in a neighborhood of the point $z_0=%
\E _\mu (\Phi )\ ,\ \Phi \in ({\cal N})_\mu ^{-1}$. Then $F^{\diamond }(\Phi
)$ defined by $S_\mu (F^{\diamond }(\Phi ))=F(S_\mu \Phi )$ exists in $(%
{\cal N})^{-1}$ .
\end{theorem}

\TeXButton{Proof}{\proof}By the characterization Theorem \ref{CharTh} $S_\mu
\Phi \in {\rm Hol}_0({\cal N}_{\Ckl })$. Then $F(S_\mu \Phi )\in {\rm Hol}_0(%
{\cal N}_{\Ckl })$ since the composition of two analytic functions is also
analytic. Again by characterization Theorem we find $F^{\diamond }(\Phi )\in
({\cal N})_\mu ^{-1}.$\TeXButton{End Proof}{\endproof}\bigskip\ 

\TeXButton{Remark }{\remark } If $F(z)=\sum_{n=0}^\infty a_k(z-z_0)^n$ then
the {\it Wick series} $\sum_{n=0}^\infty a_k(\Phi -z_0)^{\diamond n}$ (where 
$\Psi ^{\diamond n}=\Psi \diamond \ldots \diamond \Psi $ n-times converges
in $({\cal N})^{-1}$ and $F^{\diamond }(\Phi )=\sum_{n=0}^\infty a_k(\Phi
-z_0)^{\diamond n}$ holds.\bigskip\ 

\example  The above mentioned equation $\Phi \diamond X=\Psi $ can be solved
if $\E _\mu (\Phi )=S_\mu \Phi (0)\neq 0$. That implies $(S_\mu \Phi
)^{-1}\in {\rm Hol}_0({\cal N}_{\Ckl })$. Thus $\Phi ^{\diamond (-1)}=S_\mu
^{-1}\left( (S_\mu \Phi )^{-1}\right) \in ({\cal N})_\mu ^{-1}$. Then $%
X=\Phi ^{\diamond (-1)}\diamond \Psi $ is the solution in $({\cal N})_\mu
^{-1}$. For more instructive examples we refer the reader to \cite{KLS94}.

%% file: n1pos.tex
\LaTeXparent{dis3.tex}

\section{Positive distributions \label{PosDist}}

In this section we will characterize the positive distributions in $({\cal N}%
)_\mu ^{-1}$. We will prove that the positive distributions can be
represented by measures in ${\cal M}_a ({\cal N}^{\prime })$. In the case of
the Gaussian Hida distribution space $({\cal S})^{\prime }$ similar
statements can be found in works of Kondratiev \cite[b]{Ko80a} and Yokoi 
\cite{Yok90,Yok93}, see also \cite{Po87} and \cite{Lee91}. In the Gaussian
setting also the positive distributions in $({\cal N})^{-1}$ have been
discussed, see \cite{KoSW94}.

Since $({\cal N})^1={\cal E}_{\min }^1({\cal N}^{\prime })$ we say that $%
\varphi \in ({\cal N})^1\,$is positive ($\varphi \geq 0$) if and only if $%
\varphi (x)\geq 0$ for all $x\in {\cal N}^{\prime }$.

\begin{definition}
An element $\Phi \in ({\cal N})_\mu ^{-1}$ is positive if for any positive $%
\varphi \in ({\cal N})^{1\text{ }}$we have $\left\langle \!\left\langle \Phi
,\varphi \right\rangle \!\right\rangle _\mu \geq 0$ . The cone of positive
elements in $({\cal N})_\mu ^{-1}$ is denoted by $({\cal N})_{\mu ,+}^{-1}$.
\end{definition}

\begin{theorem}
\label{N1Posi} Let $\Phi \in ({\cal N})_{\mu ,+}^{-1}$ . Then there exists a
unique measure $\nu \in {\cal M}_a({\cal N}^{\prime })$ such that $\forall
\varphi \in ({\cal N})^1$ 
\begin{equation}
\label{Phiny}\left\langle \!\left\langle \Phi ,\varphi \right\rangle
\!\right\rangle _\mu =\int_{{\cal N}^{\prime }}\varphi (x)\ {\rm d}\nu (x)\ .
\end{equation}
Vice versa, any (positive) measure $\nu \in {\cal M}_a({\cal N}^{\prime })$
defines a positive distribution $\Phi \in ({\cal N})_{\mu ,+}^{-1}$ by (\ref
{Phiny}).
\end{theorem}

\TeXButton{Remarks}{\remarks } \smallskip
\\1. For a given measure $\nu $ the distribution $\Phi $ may be viewed as
the generalized Radon-Nikodym derivative $\frac{{\rm d}\nu }{{\rm d}\mu }$
of $\nu $ with respect to $\mu $. In fact if $\nu $ is absolutely continuous
with respect to $\mu $ then the usual Radon-Nikodym derivative coincides
with $\Phi .$ \smallskip
\\2. Note that the cone of positive distributions generates the same set of
measures ${\cal M}_a({\cal N}^{\prime })$ for all initial measures $\mu \in $
${\cal M}_a({\cal N}^{\prime })$. \medskip\ 

\TeXButton{Proof}{\proof}To prove the first part we define moments of a
distribution $\Phi $ and give bounds on their growth. Using this we
construct a measure $\nu $ which is uniquely defined by given moments%
\TeXButton{TeX}{\renewcommand{\thefootnote}{\fnsymbol{footnote}}}\footnote{%
Since the algebra of exponential functions is not contained in $({\cal N}%
)_\mu ^1$ we cannot use Minlos' theorem to construct the measure. This was
the method used in Yokoi's work \cite{Yok90}.}. The next step is to show
that any test functional $\varphi \in {\cal (N)}^1$ is integrable with
respect to $\nu $.

Since ${\cal P}({\cal N}^{\prime })\subset {\cal (N)}^1$ we may define
moments of a positive distribution $\Phi \in ({\cal N})_\mu ^{-1}$ by 
$$
{\rm M}_n(\xi _1,...,\xi _n)=\left\langle \!\!\!\left\langle \Phi ,\
\prod\limits_{j=1}^n\left\langle \cdot ,\xi _j\right\rangle \right\rangle
\!\!\!\right\rangle _\mu \ ,\quad \ n\in {\N},\quad \xi _j\in {\cal N,\ }%
1\leq j\leq n 
$$
$$
{\rm M}_0=\left\langle \!\left\langle \Phi ,\ \1 \right\rangle
\!\right\rangle \ \text{.} 
$$
We want to get estimates on the moments. Since $\Phi \in ({\cal H}%
_{-p})_{-q,\mu }^{-1}$ for some $p,q>0$ we may estimate as follows 
$$
\bigg|\left\langle \!\!\!\left\langle \Phi ,\left\langle x^{\otimes
n},\tbigotimes_{j=1}^n\xi _j\right\rangle \right\rangle \!\!\!\right\rangle
_\mu \bigg|\leq \left\| \Phi \right\| _{-p,-q,\mu }\left\| \left\langle
x^{\otimes n},\tbigotimes_{j=1}^n\xi _j\right\rangle \right\| _{p,q,\mu }%
\text{ .} 
$$
To proceed we use the property (\ref{(P2)}) and the estimate (\ref{MnmuNorm}%
) to obtain 
\begin{eqnarray*}
\left\| \left\langle x^{\otimes n},\tbigotimes_{j=1}^n\xi _j\right\rangle %
\right\| _{p,q,\mu }^2 %
&=& \sum_{k=0}^n\binom nk^2\left\| \left\langle P_k^\mu %
\hat \otimes {\rm M}_{n-k}^\mu ,\tbigotimes_{j=1}^n\xi _j\right\rangle %
\right\| _{p,q,\mu }^2 %
\\& \leq & \sum_{k=0}^n\binom nk^2(k!)^2 \, 2^{kq}\,|{\rm M}_{n-k}^\mu %
|_{-p}^2\prod_{j=1}^n|\xi _j|_p^2 %
\\&=& \tprod_{j=1}^n|\xi _j|_p^2\sum_{k=0}^n\binom nk^2(k!)^2\left( (n-k)! %
\right) ^2K^{2(n-k)}2^{kq} %
\\& \leq & \tprod_{j=1}^n|\xi %
_j|_p^2\;(n!)^2 \, 2^{nq}\sum_{k=0}^n2^{-(n-k)q}K^{2(n-k)} %
\\& \leq & \tprod_{j=1}^n|\xi _j|_p^2\;(n!)^2 \, 2^{nq}%
\sum_{k=0}^\infty 2^{-kq}K^{2k} 
\end{eqnarray*}
which is finite for $p,q$ large enough. Here $K$ is determined by equation (%
\ref{MnmuNorm}).

\noindent Then we arrive at 
\begin{equation}
\label{mombound}\Big|{\rm M}_n(\xi _1,...\xi _n)\Big|\leq K\ C^n\
n!\prod\limits_{j=1}^n\left| \xi _j\right| _p 
\end{equation}
for some $K,C>0$.

Due to the kernel theorem \ref{KernelTh} we then have the representation%
$$
{\rm M}_n(\xi _1,...\xi _n)=\left\langle {\rm M}^{(n)},\xi _1\otimes
...\otimes \xi _n\right\rangle \text{ ,} 
$$
where ${\rm M}^{(n)}\in ({\cal N}^{\hat \otimes n})^{\prime }$. The sequence 
$\left\{ {\rm M}^{(n)},\ n\in {\N}_0\right\} $ has the following property of
positivity: for any finite sequence of smooth kernels $\left\{ g^{(n)},n\in {%
\N}\right\} $ (i.e.,\ $g^{(n)}\in {\cal N}^{\hat \otimes n}$ and $g^{(n)}=0$%
\ $\forall \;\,n\geq n_0$ for some $n_0\in {\N}$) the following inequality
is valid 
\begin{equation}
\label{posmon}\sum_{k,j}^{n_0}\left\langle {\rm M}^{(k+j)}\;,g^{(k)}\otimes 
\overline{g^{(j)}}\right\rangle \geq 0\text{ .} 
\end{equation}
This follows from the fact that the left hand side can be written as $%
\left\langle \!\left\langle \Phi ,|\varphi |^2\right\rangle \!\right\rangle $
with%
$$
\varphi (x)=\sum_{n=0}^{n_0}\left\langle x^{\otimes n},g^{(n)}\right\rangle
,\quad x \in {\cal N}^{\prime }\text{ ,} 
$$
which is a smooth polynomial. Following \cite{BS71,BeKo88} inequalities (\ref
{mombound}) and (\ref{posmon}) are sufficient to ensure the existence of a
uniquely defined measure $\nu $ on $({\cal N}^{\prime },{\cal C}_\sigma (%
{\cal N}^{\prime }))$, such that for any $\varphi \in {\cal P}({\cal N}%
^{\prime })$ we have 
$$
\left\langle \!\left\langle \Phi ,\varphi \right\rangle \!\right\rangle _\mu
=\int_{{\cal N}^{\prime }}\varphi (x)\ {\rm d}\nu (x)\text{ .} 
$$

From estimate (\ref{mombound}) we know that $\nu \in {\cal M}_a({\cal N}%
^{\prime })$. Then Lemma \ref{equiLemma} shows that there exists $%
\varepsilon >0\,,\,p\in \N $ such that $\exp (\varepsilon |x|_{-p})$ is $\nu
-$integrable. Corollary \ref{phi(z)Betrag} then implies that each $\varphi
\in ({\cal N})^1$ is $\nu $-integrable.\bigskip\ 

Conversely let $\nu \in {\cal M}_a({\cal N}^{\prime })$ be given. Then the
same argument shows that each $\varphi \in ({\cal N})^1$ is $\nu $%
-integrable and from Corollary \ref{phi(z)Betrag} we know that%
$$
\left| \int_{{\cal N}^{\prime }}\varphi (x){\rm d}\nu (x)\right| \leq
\,C\,\left\| \varphi \right\| _{p,q,\mu }\int_{{\cal N}^{\prime }}\exp
(\varepsilon |x|_{-p})\,{\rm d}\nu (x) 
$$
for some $p,q\in \N \,,\,C>0$. Thus the continuity of $\varphi \mapsto \int_{%
{\cal N}^{\prime }}\varphi \;{\rm d}\nu $ is established, showing that $\Phi 
$ defined by equation (\ref{Phiny}) is in $({\cal N})_{\mu ,+}^{-1}$.%
\TeXButton{End Proof}{\endproof}

%% file: change.tex
\LaTeXparent{nga4.tex}

\section{Change of measure}

Suppose we are given two measures $\mu ,\hat \mu \in {\cal M}_a({\cal N}%
^{\prime })$ both satisfying Assumption 2. Let a distribution $\hat \Phi \in
({\cal N})_{\hat \mu }^{-1}$ be given. Since the test function space $({\cal %
N})^1$ is invariant under changes of measures in view of Theorem \ref
{N1E1min}, the continuous mapping 
$$
\varphi \mapsto \langle \!\langle \hat \Phi ,\varphi \rangle \!\rangle
_{\hat \mu }\ ,\qquad \varphi \in ({\cal N})^1 
$$
can also be represented as a distribution $\Phi \in ({\cal N})_\mu ^{-1}$.
So we have the implicit relation $\Phi \in ({\cal N})_\mu
^{-1}\leftrightarrow \hat \Phi \in ({\cal N})_{\hat \mu }^{-1}$ defined by 
$$
\langle \!\langle \hat \Phi ,\varphi \rangle \!\rangle _{\hat \mu }=\langle
\!\langle \Phi ,\varphi \rangle \!\rangle _\mu \ . 
$$
This section will provide formulae which make this relation more explicit in
terms of re-decomposition of the $\Q $-series. First we need an explicit
relation of the corresponding $\p $-systems.

\begin{lemma}
Let $\mu ,\hat \mu \in {\cal M}_a({\cal N}^{\prime })$ then 
$$
P_n^\mu (x)=\sum_{k+l+m=n}\frac{n!}{k!\,l!\,m!}P_k^{\hat \mu }(x)\hat
\otimes P_l^\mu (0)\hat \otimes {\rm M}_m^\mu \ . 
$$
\end{lemma}

\TeXButton{Proof}{\proof}Expanding each factor in the formula%
$$
e_\mu (\theta ,x)=e_{\hat \mu }(\theta ,x)l_\mu ^{-1}(\theta )l_{\hat \mu
}(\theta )\ , 
$$
we obtain%
$$
\sum_{n=0}^\infty \frac 1{n!}\langle P_n^\mu (x),\theta ^{\otimes n}\rangle
=\sum_{k,l,m=0}^\infty \frac 1{k!\,l!\,m!}\langle P_k^\mu (x)\otimes
P_l^{\hat \mu }(0)\otimes {\rm M}_m^\mu ,\theta ^{\otimes (k+l+m)}\rangle \
. 
$$
A comparison of coefficients gives the above result.\TeXButton{End Proof}
{\endproof}\bigskip\ 

An immediate consequence is the next reordering lemma.

\begin{lemma}
Let $\varphi \in ({\cal N})^1$ be given. Then $\varphi $ has representations
in $\p  ^\mu $-series as well as $\p  ^{\hat \mu }$-series:%
$$
\varphi =\sum_{n=0}^\infty \langle P_n^\mu ,\varphi ^{(n)}\rangle
=\sum_{n=0}^\infty \langle P_n^{\hat \mu },\hat \varphi ^{(n)}\rangle  
$$
where $\varphi ^{(n)},\hat \varphi ^{(n)}$ $\in {\cal N}_{\Ckl }^{\hat
\otimes n}$ for all $n\in \N _0$, and the following formula holds: 
\begin{equation}
\label{phiHut}\hat \varphi ^{(n)}=\sum_{l,m=0}^\infty \frac{(l+m+n)!}{%
l!\,m!\,n!}\left( P_l^\mu (0)\hat \otimes {\rm M}_m^{\hat \mu },\varphi
^{(l+m+n)}\right) _{{\cal H}^{\otimes (l+m)}}\ .
\end{equation}
\end{lemma}

\noindent Now we may prove the announced theorem.

\begin{theorem}
Let $\hat \Phi =\sum_{n=0}^\infty \langle Q_n^{\hat \mu },\hat \Phi
^{(n)}\rangle \in ({\cal N)}_{\hat \mu }^{-1}$. Then $\Phi
=\sum_{n=0}^\infty \langle Q_n^\mu ,\Phi ^{(n)}\rangle $ defined by 
$$
\langle \!\langle \Phi ,\varphi \rangle \!\rangle _\mu =\langle \!\langle
\hat \Phi ,\varphi \rangle \!\rangle _{\hat \mu }\ ,\qquad \varphi \in (%
{\cal N})^1 
$$
is in $({\cal N})_\mu ^{-1}$ and the following relation holds%
$$
\Phi ^{(n)}=\sum_{k+l+m=n}\frac 1{l!\,m!}\hat \Phi ^{(k)}\hat \otimes
P_l^\mu (0)\hat \otimes {\rm M}_m^{\hat \mu } 
$$
\end{theorem}

\TeXButton{Proof}{\proof}We can insert formula (\ref{phiHut}) in the formula 
$$
\sum_{n=0}^\infty n!\, \langle \Phi ^{(n)},\varphi ^{(n)}\rangle
=\sum_{n=0}^\infty n!\, \langle \hat \Phi ^{(n)},\hat \varphi ^{(n)}\rangle 
$$
and compare coefficients again.\TeXButton{End Proof}{\endproof}

%% file: NGA5.bbl
\begin{thebibliography}{KLPSW94}
\addcontentsline{toc}{section}{References}
\addtolength{\itemsep}{-1mm}

\bibitem[AKS93]{AKS93}  Albeverio, S., Kondratiev, Yu.G. and Streit, L.
(1993), {\it How to generalize White Noise Analysis to Non-Gaussian Spaces}.
In: 'Dynamics of Complex and Irregular Systems'. Eds.: Ph. Blanchard et al.,
World Scientific.

\bibitem[ADKS94]{ADKS94}  Albeverio, S., Daletzky, Y., Kondratiev, Yu. G. and
Streit, L. (1994), {\it Non-Gaussian infinite dimensional analysis, }%
preprint, to appear in J. Func. Anal..

\bibitem[BeS95]{BeS95}  Benth, F. and Streit, L. (1995), {\it The Burgers
Equation with a Non-Gaussian Random Force. }UMa preprint.

\bibitem[BeKo88]{BeKo88}  Berezansky, Yu. M. and Kondratiev, Yu. G. (1988), 
{\it Spectral Methods in Infinite-Dimensional Analysis}, (in Russian),
Naukova Dumka, Kiev. English translation, 1995, Kluwer Academic Publishers,
Dordrecht.

\bibitem[BeLy93]{BeLy93}  Berezansky, Yu.M. and Lytvynov, E.V. (1993), {\it %
Generalized White Noise Analysis connected with pertubed field operators},
Dopovidy AN Ukrainy, No {\bf 10}.

\bibitem[BS71]{BS71}  Berezansky, Yu. M. and Shifrin, S.N. (1971), {\it The
generalized degree symmetric Moment Problem,} Ukrainian Math. J. 23 N3,
247-258.

\bibitem[Bo76]{Bo76}  Bourbaki, N. (1976), {\it Elements of mathematics.
Functions of a real variable. }Addison-Wesley.

\bibitem[Da91]{Da91}  Daletsky, Yu.L. (1991),{\it \ A biorthogonal analogy
of the Hermite polynomials and the inversion of the Fourier transform with
respect to a non Gaussian measure}, Funct. Anal. Appl. {\bf 25}, 68-70.

\bibitem[Di81]{Di81}  Dineen, S. (1981), {\it Complex Analysis in Locally
Convex Spaces,} Mathematical Studies 57, North Holland, Amsterdam.

\bibitem[GV68]{GV68}  Gel'fand, I. M. and Vilenkin, N.Ya. (1968),{\it \
Generalized Functions}, Vol. IV, Academic Press, New York and London.

\bibitem[Hi75]{Hi75}  Hida, T. (1975), {\it Analysis of Brownian
Functionals, }Carleton Math. Lecture Notes No. 13, Carleton.

\bibitem[Hi80]{Hi80}  Hida, T. (1980), {\it Brownian Motion}. Springer, New
York.

\bibitem[HKPS93]{HKPS93}  Hida, T., Kuo, H.H., Potthoff, J. and Streit, L.
(1993),{\it \ White Noise. An infinite dimensional calculus}. Kluwer,
Dordrecht.

\bibitem[Ito88]{Ito88}  Ito, Y. (1988), {\it Generalized Poisson Functionals.%
} Prob. Th. Rel. Fields {\bf 77} 1-28.

\bibitem[IK88]{IK88}  Ito, Y. and Kubo, I. (1988), {\it Calculus on Gaussian and Poisson
White Noises.} Nagoya Math. J. {\bf 111} 41-84.

\bibitem[Ko78]{Ko78}  Kondratiev, Yu.G. (1978), {\it Generalized functions
in problems of infinite dimensional analysis}. Ph.D. thesis, Kiev University.

\bibitem[Ko80a]{Ko80a}  Kondratiev, Yu.G. (1980), {\it Spaces of entire
functions of an infinite number of variables, connected with the rigging of
a Fock space. }In: 'Spectral Analysis of Differential Operators.' Math.
Inst. Acad. Sci. Ukrainian SSR, p. 18-37. English translation: Selecta Math.
Sovietica {\bf 10} (1991), 165-180.

\bibitem[Ko80b]{Ko80b}  Kondratiev, Yu.G. (1980), {\it Nuclear spaces of
entire functions in problems of infinite dimensional analysis.} Soviet Math.
Dokl. {\bf 22}, 588-592.

\bibitem[KLPSW94]{KLPSW94}  Kondratiev, Yu.G., Leukert, P., Potthoff, J.,
Streit, L., Westerkamp, W. (1994), {\it Generalized Functionals in Gaussian
Spaces -- the Characterization Theorem Revisited. }Manuskripte 175/94, Uni
Mannheim.

\bibitem[KLS94]{KLS94}  Kondratiev, Yu.G., Leukert, P., Streit, L. (1994),%
{\it \ Wick Calculus in Gaussian Analysis,} BiBoS preprint 637, to appear in
Acta Applicandae Mathematicae.

\bibitem[KoSa78]{KoSa78}  Kondratiev, Yu.G. and Samoilenko, Yu.S. (1978), 
{\it Spaces of trial and generalized functions of an infinite number of
variables}, Rep. Math. Phys. {\bf 14}, No.3, 325-350.

\bibitem[KoS93]{KoS92}  Kondratiev, Yu.G. and Streit, L. (1993), {\it Spaces
of White Noise distributions: Constructions, Descriptions, Applications.} I.
Rep. Math. Phys. {\bf 33}, 341-366.

\bibitem[KoSW95]{KoSW94}  Kondratiev, Yu.G. , Streit, L. and Westerkamp,
W.(1995), {\it A Note on Positive Distributions in Gaussian Analysis},
Ukrainian Math. J. {\bf 47} No. 5.

\bibitem[KoTs91]{KoTs91}  Kondratiev, Yu.G. and Tsykalenko T.V. (1991), {\it %
Dirichlet Operators and Associated Differential Equations.} Selecta Math.
Sovietica {\bf 10}, 345-397.

\bibitem[KMP65]{KMP65}  Kristensen, P., Mejlbo, L., and Poulsen, E.T.
(1965), {\it Tempered Distributions in Infinitely Many Dimensions. I.
Canonical Field Operators. }Commun. math. Phys. {\bf 1}, 175--214.

\bibitem[Kuo92]{Kuo92}  Kuo, H.-H. (1992), {\it Lectures on white noise
analysis}. Soochow J. Math. {\bf 18}, 229-300.

\bibitem[Lee91]{Lee91}  Lee, Y.-J. (1991), {\it Analytic Version of Test
Functionals, Fourier Transform and a Characterization of Measures in White
Noise Calculus. }J. Funct. Anal. {\bf 100}, 359-380.

\bibitem[Lu70]{Lu70}  Lukacs, E. (1970), {\it Characteristic Functions}, 2nd
edition, Griffin, London.

\bibitem[\O k95]{Ok94}  \O ksendal, B. (1995),{\it \ Stochastic Partial
Differential Equations and Applications to Hydrodynamics. }In:
`Stochastic Analysis and Applications to Physics' Eds.: A.I. Cardoso
et al; Kluwer, Dordrecht.

\bibitem[Ou91]{Ou91}  Ouerdiane, H. (1991),{\it \ Application des m\'ethodes
d'holomorphie et de distributions en dimension quelconque \'a l'analyse sur
les espaces Gaussiens}. BiBoS preprint 491.

\bibitem[Po87]{Po87}  Potthoff, J. (1987), {\it On positive generalized
functionals. }J. Funct. Anal. {\bf 74}, 81-95.

\bibitem[PS91]{PS91}  Potthoff, J. and Streit, L. (1991), {\it A
characterization of Hida distributions.} J. Funct. Anal. {\bf 101}, 212-229.

\bibitem[Sch71]{Sch71}  Schaefer, H.H. (1971), {\it Topological Vector
Spaces, }Springer, New York.

\bibitem[Sk74]{Sk74}  Skorohod, A.V. (1974), {\it Integration in Hilbert
Space, }Springer, Berlin.

\bibitem[Va95]{Va95}  V\aa ge, G. (1995), {\it Stochastic Differential
Equations and Kondratiev Spaces. }Ph.D. thesis, Trondheim University. 

\bibitem[VGG75]{VGG75}  Vershik, A.M., Gelfand, I.M. and Graev, M.I. (1975), 
{\it Representations of diffeomorphisms groups}. Russian Math. Surveys {\bf %
30}, No 6, 3-50.

\bibitem[Yok90]{Yok90}  Yokoi, Y.(1990), {\it Positive generalized white
noise functionals}. Hiroshima Math. J. {\bf 20}, 137-157.

\bibitem[Yok93]{Yok93}  Yokoi, Y. (1993), {\it Simple setting for white
noise calculus using Bargmann space and Gauss transform.} Preprint.

\end{thebibliography}
